\documentclass[12pt]{article}
\usepackage{amssymb,epsfig,amsfonts,epic,graphics,pictex} 
\begin{document}

\newtheorem{lem}{Lemma}[section]
\newtheorem{prop}{Proposition}
\newtheorem{con}{Construction}[section]
\newtheorem{defi}{Definition}[section]
\newtheorem{coro}{Corollary}[section]
\newcommand{\hf}{\hat{f}}
\newtheorem{fact}{Fact}[section]
\newtheorem{theo}{Theorem}
\newcommand{\Br}{\Poin}
\newcommand{\Cr}{{\bf Cr}}
\newcommand{\dist}{{\rm dist}}
\newcommand{\diam}{\mbox{diam}\, }
\newcommand{\mod}{{\rm mod}\,}
\newcommand{\compose}{\circ}
\newcommand{\dbar}{\bar{\partial}}
\newcommand{\Def}[1]{{{\em #1}}}
\newcommand{\dx}[1]{\frac{\partial #1}{\partial x}}
\newcommand{\dy}[1]{\frac{\partial #1}{\partial y}}
\newcommand{\Res}[2]{{#1}\raisebox{-.4ex}{$\left|\,_{#2}\right.$}}
\newcommand{\sgn}{{\rm sgn}}

\newcommand{\CC}{\mathbb{C}}
\newcommand{\D}{{\bf D}}
\newcommand{\Dm}{{\bf D_-}}
\newcommand{\RR}{\mathbb{R}}
\newcommand{\NN}{\mathbb{N}}
\newcommand{\HH}{\mathbb{H}}
\newcommand{\Z}{{\bf Z}}
\newcommand{\tr}{\mbox{Tr}\,}
\newcommand{\R}{{\bf R}}
\newcommand{\C}{{\bf C}}

\newenvironment{nproof}[1]{\trivlist\item[\hskip \labelsep{\bf Proof{#1}.}]}
{\begin{flushright} $\square$\end{flushright}\endtrivlist}
\newenvironment{proof}{\begin{nproof}{}}{\end{nproof}}

\newenvironment{block}[1]{\trivlist\item[\hskip \labelsep{{#1}.}]}{\endtrivlist}
\newenvironment{definition}{\begin{block}{\bf Definition}}{\end{block}}

\newtheorem{conjec}{Conjecture}

\newtheorem{com}{Comment}
\font\mathfonta=msam10 at 11pt
\font\mathfontb=msbm10 at 11pt
\def\Bbb#1{\mbox{\mathfontb #1}}
\def\lesssim{\mbox{\mathfonta.}}
\def\suppset{\mbox{\mathfonta{c}}}
\def\subbset{\mbox{\mathfonta{b}}}
\def\grtsim{\mbox{\mathfonta\&}}
\def\gtrsim{\mbox{\mathfonta\&}}

\newcommand{\Poin}{{\bf Poin}}
\newcommand{\Bo}{\Box^{n}_{i}}
\newcommand{\Di}{{\cal D}}
\newcommand{\gd}{{\underline \gamma}}
\newcommand{\gu}{{\underline g }}
\newcommand{\ce}{\mbox{III}}
\newcommand{\be}{\mbox{II}}
\newcommand{\F}{\cal{F}}
\newcommand{\Ci}{\bf{C}}
\newcommand{\ai}{\mbox{I}}
\newcommand{\dupap}{\partial^{+}}
\newcommand{\dm}{\partial^{-}}
\newenvironment{note}{\begin{sc}{\bf Note}}{\end{sc}}
\newenvironment{notes}{\begin{sc}{\bf Notes}\ \par\begin{enumerate}}%
{\end{enumerate}\end{sc}}
\newenvironment{sol}
{{\bf Solution:}\newline}{\begin{flushright}
{\bf QED}\end{flushright}}

\title{Dynamics and Universality of Unimodal Mappings with Infinite 
Criticality}

\author{Genadi Levin
\thanks{Both authors were supported by Grant No. 2002062
from the United States-Israel Binational Science Foundation
(BSF), Jerusalem, Israel.}\\
\small{Dept.\ of Math.}\\
\small{Hebrew University}\\
\small{Jerusalem 91904, ISRAEL}\\
\small{\tt levin@math.huji.ac.il}\\
\and
Grzegorz \'{S}wia\c\negthinspace tek
\thanks{Partially supported by NSF grant DMS-0245358.}\\
\small{Dept.\ of Math.}\\
\small{Penn State University}\\
\small{University Park, PA 16802, USA}\\
\small{\tt swiatek@math.psu.edu}}
\normalsize
\maketitle
\abstract{We consider infinitely renormalizable unimodal mappings with
topological type which is periodic under renormalization. We study the limiting
behavior of fixed points of the renormalization operator 
as the order of the critical point increases to infinity. It
is shown that a limiting dynamics exists, with
a critical point that is flat, but still having a well-behaved
analytic continuation to a neighborhood of the real interval pinched
at the critical point. We study the dynamics of limiting maps and
prove their rigidity. In particular, the sequence of fixed points of
renormalization for finite criticalities converges, uniformly on the
real domain, to a mapping of the limiting type.}

\section{Introduction}
\subsection{Overview of the problem.}
Universality for unimodal mappings 
was discovered by Feigenbaum~\cite{feig0},~\cite{feig1} and 
Coullet-Tresser~\cite{CT}
in the case of period doubling, initially purely on the basis of
numerical observation. For our purposes, the problem can be stated as
follows. We consider  
mappings $H :\: [0,1] \rightarrow
[0,1]$ in the form 
\[ H(x) = |E(x)|^{\ell} \]
where $\ell>1$ is a real number and 
$E$ is a smooth mapping with strictly negative
derivative on $[0,1]$ which maps $0$ to $1$ and $1$ to a point inside
$(-1,0)$. Then $H$ is unimodal with the minimum at some 
$x_{0} = E^{-1}(0) \in (0,1)$ and $x_0$ is the critical
point of order $\ell$. The celebrated {\em Feigenbaum functional equation}
is 
\begin{equation}\label{equ:14fa,1}
\tau H^2(x) = H(\tau x) 
\end{equation}
for $x \in [0,\tau^{-1}]$. The equation needs to be solved for $H$ and
then necessarily $\tau^{-1} = H(1) = H^2(0)$.    

The original discovery was that the the solution to Feigenbaum's
functional equation can be found by iterating the following {\em
renormalization operator} 
\[ {\cal R}(H)(x) :=  \frac{H^2(H(1)\cdot x)}{H(1)} \]
which can be seen as a step in the method of successive approximations
for solving Equation~(\ref{equ:14fa,1}).  Note that ${\cal
R}(H)$ satisfies conditions imposed in the preceding paragraph
provided that $H(1) < x_0$ and then $\cal R$ can be applied again. Universality
means that as soon as
${\cal R}^n(H)$ remains in the class described above for all $n$, this sequence
will converge to a limit $H_{\ell}$ also in the same class, and
$H_{\ell}$  is a
solution to Feigenbaum's functional equation. Moreover, this limit is independent of the
initial guess $H$, except for the rank of criticality $\ell$.  

The early thrust of the theory was toward actually solving
Feigenbaum's equation and finding constants $\tau_{\ell}$ for small
values of $\ell$. Next, rigorous 
computer-assisted proofs were developed,
see~\cite{lanf0},~\cite{lanf},~\cite{coep}. Later, the problem was
generalized to include versions of Equation~(\ref{equ:14fa,1}) which
involve a higher iterate of $H$ replacing the second one. 

At that point the need for a more theoretical approach to the problem became
obvious. First, one could not re-run computer estimates in all
infinitely many cases to which the theory seemed to apply. Secondly, 
while computer-assisted proofs showed the emergence of universal
constants and functions, it still did not explain 
qualitative reasons of the phenomena. The program for solving
renormalization conjectures purely with tools of dynamical systems
theory was formulated by D. Sullivan in the
mid-1980s. Its salient feature was strong reliance on {\em complex}
dynamics of analytic continuations of real maps.  
This approach took some time to develop, but has been highly
successful in the end, 
see~\cite{Su},~\cite{profesorus}, ~\cite{wisnia}, ~\cite{mcmdoklad}. 
In particular, for each $\ell$ which is an even integer the 
existence of a solution $H_{\ell}$ to
Equation~(\ref{equ:14fa,1}), unique for the order of criticality
$\ell$,  has been rigorously established. 

This paper is concerned with the case when $\ell$ increases to
$\infty$. Originally, interest in this problem came from mathematical
physics literature, see~\cite{EW},~\cite{wiele},~\cite{thompson},
~\cite{briggs}.
One motivation came from the expectation that the problem could shed
light on other, more complicated,  limit problems of
statistical and quantum physics. Another reason was the obvious computational
challenge of working with Equation~(\ref{equ:14fa,1}) for large
$\ell$. With such $\ell$, the renormalization operator cannot be
iterated for very long 
because of finite accuracy and hence different procedures were needed
for solving the equation. Papers quoted here all successfully dealt
with this challenge obtaining consistent estimates for
$\lim_{\ell\rightarrow \infty} \tau_{\ell} \approx 30$, for example. 
Their methods were cast in varying language, but were all based on the
fact that functions $H_{\ell}$ for $\ell \rightarrow \infty$ approach
the Fatou coordinate of a certain parabolic fixed point. In addition
to developing a numerical approach, paper~\cite{EW} contained a
rigorous computer-assisted proof of the existence of a limiting
function $H_{\infty}$ which solved Equation~(\ref{equ:14fa,1}) and was the
limit of fixed-point transformations $H_\ell$ for finite $\ell$,
\[ H_{\ell}(x) = |E_{\ell}(x)|^{\ell}\; . \]
It is actually curious that such a
limit may exist at all. Here, it happens because $E_{\ell}(x)$ for some
fixed $x\neq x_0$ will tend to $-1$ or $1$ at a rate proportional to
$\ell^{-1}$. 

The second source of interest was the study of metric attractors of
real and complex maps. In~\cite{bruke}, the first example of an exotic
attractor was shown for a unimodal map. The key estimate of the paper
was obtained by adjusting the order of the critical point to a
sufficiently high value. That work was followed by a program of
S. Van Strien and T. Nowicki for showing the existence of a similar
attractor for a complex polynomial, which would imply that the Julia
set of such a polynomial has positive measure. While that program has
not been followed to a successful completion, partial progress was based
on choosing sufficiently high criticality and studying limits when it
tended to $\infty$, see~\cite{NStr}.

\paragraph{Contribution of this paper.}
We provide an analytic, not computer-assisted, proof of the existence
and uniqueness 
of the solution to the generalized Equation~(\ref{equ:14fa,1}) in a
so-called EW-class
of mappings with an infinitely flat critical point, see next subsection. 
Maps from the EW-class cover all topological equivalence classes which
contain infinitely renormalizable transformations from the quadratic
family and are periodic under the renormalization.
Since our class contains limits of sequences
$H_{\ell}$ as $\ell \rightarrow \infty$ of fixed points of
renormalization for finite $\ell$, the uniqueness means that
$H_{\ell}$ actually converge to that limiting infinitely flat
dynamics. 
Ultimately the EW-class is precisely the class  of the limiting
maps: for every
order type $\aleph$, as defined in detail further, 
the EW-class contains one and only one map
$H_\aleph$ of this type, and, moreover, $H_\aleph$ is the limit
of any sequence of fixed-point maps 
$H_{\ell}$ of the same type $\aleph$,
as the criticality $\ell$ tends to infinity along
real numbers.

We also study basic properties of the limiting map as
complex dynamical system. 

Technically, our approach is based on the rigidity of towers in the
sense of~\cite{profesorus}. The class of complex maps we
are working with is quite different from polynomial-like mappings
studied for finite 
even integer $\ell$. 
Moreover, the sequence of maps $H_{\ell}$
is not generated by any identifiable operator in a functional space.
In spite of these significant differences with the standard setting, 
the basic approach still works. It looks
likely that it should also work for other types of dynamics such as
circle homeomorphisms or Fibonacci induced maps.

Main results of the paper are contained  in 
Theorems~\ref{theo:13up,1}-~\ref{theo:13up,2}.

\subsection{Statement of main results.}
We will say that two finite sequences $(u_i)_{i=1}^p$ and
$(u'_i)_{i=1}^{p'}$ have the same {\em order type} provided that $p=p'$ and 
$u_i < u_j$ iff $\varepsilon u'_i < \varepsilon u'_j$ for all 
$i,j = 1,\cdots, p$ and a fixed constant $\varepsilon$.  The order
type is an equivalence class of this relation, typically denoted with
a Hebrew letter, and then $|\aleph|$ will mean the length of a
sequence in $\aleph$. We will consider infinitely renormalizable maps
with periodic combinatorics given by some order type $\aleph$. This
means that for every $n$ there is a restrictive interval of period
$|\aleph|^n$ and the  order type of points $x_0,
f^{|\aleph|^{n-1}}(x_0),\cdots,  f^{(|\aleph|-1)|\aleph|^{n-1}}(x_0)$ is
$\aleph$. Here,  $x_0$ is the critical point of the univalent map
$f$. 

Unimodal maps will be denoted by $H$, often with a subscript indicating the
order of the critical point. 
That is, $H_\ell$ is assumed to be in the following form:
$H_{\ell}(x)=|E_\ell(x)|^\ell$, where $E_\ell :\: [0,1] \rightarrow
\RR$ is a $C^2$-diffeomorphism onto its image. 
Unimodal maps are normalized so that $H([0,1])
= [0,1]$,  $H(0)=1$ and the global strict minimum $0$ is attained in
$(0,1)$.   They are further assumed to be infinitely renormalizable
with some combinatorial order type $\aleph$ and to satisfy the fixed
point equation:

\begin{equation}\label{equ:10fa,2}
 \tau H^{|\aleph|}(x) = H(\tau x) \; .
\end{equation}
with $\tau>0$. 
By renormalization theory, see~\cite{Su}, a fixed point $H_{\ell}$ for
any $\ell>1$ can be represented as $|E_{\ell}|^\ell$ with $E_{\ell}$ which is
a diffeomorphism in the Epstein class:
\begin{defi}\label{defi:eps}
A diffeomorphism $E$ of a real interval $T'$ onto another real
interval $T$ is said to be in the Epstein class if
the inverse map $E^{-1}: T\to T'$ extends to a univalent
map $E^{-1} :\: (\CC \setminus \RR) \cup T' \rightarrow
(\CC\setminus\RR) \cup T$.
\end{defi}
Our first main result is the following.
\begin{theo}\label{theo:13up,1}
Let us fix an order type $\aleph$ and consider a sequence
$H_{\ell_m}$, with $\ell_m$ real, 
of unimodal maps which are infinitely renormalizable with
periodic combinatorics of type $\aleph$ and satisfy the fixed
point equation~(\ref{equ:10fa,2}), each with its own scaling constant
$\tau_{m} > 1$.  

If $\lim_{m\rightarrow\infty} \ell_m = \infty$, then $H_{\ell_m}$ converge as  
$m\rightarrow \infty$, uniformly on $[0,1]$, to a unimodal function
$H$. Also, $\lim_{m\rightarrow \infty} \tau_{m} = \tau>1$ exists,
and $H, \tau$ satisfy the fixed point equation~(\ref{equ:10fa,2}). 
\end{theo}

\paragraph{Eckmann-Wittwer class.}
One can say more about the analytic continuation of $H$. Not
only does the analytic continuation provide more information about the
limit, but is also crucial for our proof which relies on holomorphic
dynamics. Different from the theory for finite $\ell$ 
in which the analytic continuations
of limits belong to the well-known class of polynomial-like mappings, 
$H$ belongs to a limiting class of mappings with a flat critical point. 
\begin{defi}\label{defi:10fa,1}
Let $H$ be a smooth unimodal map defined from the interval $[0,1]$
into itself, with the minimum at some point $x_0 \in (0,1)$. Suppose
that it is normalized so that $H(x_0) = 0, H(0) = 1$ and the orbit
$x_0, \cdots, H^{p-1}(0)$ has order type $\aleph$. Then we will say
that $H$ belongs to the {\em Eckmann-Wittwer} class, EW-class for
short, with combinatorial type $\aleph$, provided that the following 
conditions hold.

\begin{enumerate}    
\item
$\tau H^{p} (x) = H(\tau x)$
for some scaling constant $\tau > 1$ and every $0 \leq x \leq \tau^{-1}$. 
\item
$H$ has analytic continuation to the union of two topological disks 
$U_{-}$ and $U_{+}$ and this analytic continuation will also be
denoted with $H$. 
\item
For some $R>1$, $H$ restricted to either $U_+$ or $U_-$ is a covering
(unbranched) of the punctured disk $V := D(0, R) \setminus \{0\}$ and
$\overline{U_+ \cup U_-} \subset D(0,R)$.  
\item
\[ H(z) = \lim_{m\rightarrow\infty} \left((E_{\ell_m}(z))^2\right)^{\ell_m/2}\; ,\]
where $\ell_m \rightarrow\infty$, for each $m$ the map $E_{\ell_m}$ is
a diffeomorphism in the Epstein
class, normalized so that $E_{\ell_m}(0) = 1$ and $E_{\ell_m}(1) \in
(-1,0)$.  It is understood that $w^{\ell_m/2}$ is the principal
branch defined on the plane slit along the negative half-line and that
for   every compact subset $K$ of $U_+ \cup U_-$, the right hand side
of the equality is well defined on  $K$ for almost all $m$ 
with uniform convergence on $K$.    
\item
$U_-$ contains the interval $[ b_0', x_0)$ and $U_+$ contains the
interval $(x_0, b_0]$ where $b_0'<0$, $1<b_0<R$, $H(b_0) = H(b_0') =
b_0$ and $H'(b_0) > 1$.  
\item
$U_{\pm}$ are both symmetric with respect to the real axis and their
closures intersect exactly at $x_0$.  
\item
The mapping $G(x) := H^{p-1}(\tau^{-1} x)$ fixes $x_0$ 
and $G^2$ has the following power
series expansion at $x_0$:
\[ G^2(x) = x - \epsilon (x-x_0)^3 + O(|x-x_0|^4) \]
with $\epsilon> 0$.  
\end{enumerate}

\end{defi}

\begin{theo}\label{theo:13up,2}
For every sequence $H_{\ell_m}$ as described in the hypothesis of
Theorem~\ref{theo:13up,1},  the limiting function $H$ 
belongs to the Eckmann-Wittwer class. 
\end{theo} 

The dynamics of maps in the EW-class is studied in this paper starting
from Section~\ref{sec:15fa,1}. 

In particular, we introduce the Julia set of EW-class maps.

\

Our last result is a straightening theorem for the EW-class.

As it follows from the Straightening Theorem for polynomial-like
maps [DH], any map $H_{\ell, \aleph}$, if $\ell$ is an even integer, is
quasi-conformally conjugate  to a polynomial 
$z\mapsto z^\ell+c_{\ell, \aleph}$
in neighborhoods of their Julia sets. Here we prove
that limit maps $H_{\aleph}$ are quasi-conformally conjugate
to maps of the form $f(z)=\exp(-c(z-a)^{-2})$.

\begin{theo}\label{s}
For every map $H:U_-\cup U_+\to V$ of the EW-class 
there exists a map of the form
$f(z)=\exp(-c(z-a)^{-2})$ with some real $a,c>0$, such that $H$ and
$f$ are hybrid equivalent, i.e.
there exists a quasi-conformal homeomorphism of the
plane $h$, such that
$$h\circ H=f\circ h$$
on $U_-\cup U_+$ and $\partial h/\partial \bar z=0$ a.e. on the Julia set of
$H$. Moreover, $h$ maps the Julia set of $H$ onto the Julia set of $f$.
\end{theo} 

See last Section for the proof and comments.

\subsection{Plan of the proof.}
Theorems~\ref{theo:13up,1} and~\ref{theo:13up,2} follow immediately
from the following two statements. 

\begin{theo}\label{theo:10fa,1}
Consider a sequence of fixed-point maps $H_{\ell_m}$ with scaling
constants $\tau_m$, all of combinatorial type $\aleph$ and
satisfying the hypothesis of
Theorem~\ref{theo:13up,1}. Let $x_m$
denote the critical point of $H_{\ell_m}$.

Then, there is a subsequence $m_p$ such that $x_{m_p} \rightarrow x_0$,
$\tau_{m_p} \rightarrow \tau$ and $H_{\ell_{m_p}} \rightarrow H$, 
where $H$ belongs to  the EW-class with combinatorial type $\aleph$, 
critical point at $x_0$ and the scaling constant $\tau$. The
convergence to $H$ is uniform on the interval $[0,1]$.  
\end{theo}

\begin{theo}\label{theo:10fa,2}
Let $H_1$ and $H_2$ be two maps belonging to the EW-class with the same
combinatorial type $\aleph$. Then $H_1$ = $H_2$. 
\end{theo}

Theorem~\ref{theo:10fa,1} follows from compactness of the family
$\{H_{\ell_m}\}$, which in turn follows from
real and complex bounds.
Further examination of limit maps shows that they belong to the EW-class.

To prove Theorem~\ref{theo:10fa,2} we follow the strategy
of~\cite{Su} as realized in~\cite{profesorus},  
despite of the fact that all the basic
``starting conditions'' of this approach break down
in a transparent way for limit maps in the EW-class. 
For example, if $H$ belongs to the EW-class, then:
\begin{itemize}
\item 
as a real map, $H$ has a {\it flat} critical point (as we will
presently argue) and many techniques do not apply, not even a 
 ``no wandering interval theorem'' can be taken for granted;
\item 
most strikingly, in spite of bounded
combinatorics, the geometry of the postcritical set
of $H$ is not bounded, and, therefore, 
known methods of constructing quasi-conformal
conjugacies do not work;
\item 
as a complex map, $H$ is {\it not} extended 
holomorphically through a neighborhood of its critical point; 
in particular, neither Fatou-Julia-Baker theory for meromorphic maps 
nor Sullivan-Douady-Hubbard
theory~\cite{sulcio},~\cite{DH} of polynomial-like maps is applicable.
\end{itemize}
 
Nevertheless, the proof~\cite{profesorus} can be adapted. We consider
a {\it tower}
generated by $H$, prove that it has needed chaotic properties, 
and derive the rigidity of the tower by showing that it cannot support
an invariant line-field.

In the sequel, the combinatorics $\aleph$ is fixed, and
we omit sometimes the index $\aleph$. Also,  $p$ will be used to
denote the cardinality of $\aleph$.

\paragraph{A further comment on the EW-class.}
EW-class plays a role in the proof which is somewhat analogous to the
impact of polynomial-like mappings in the standard theory. Both
classes share a fundamental ``expansion'' characteristic: namely a
smaller domain provides a covering of a larger one with the critical
value removed. However, the critical point in the EW-class is not in
the domain of analyticity. 

Assume now that $H$ belongs to the EW-class.
By the functional equation~(\ref{equ:10fa,2}), 
\[ \tau^{-1} H(z) = H(G(z)) \]
 which initially
holds for $z\in [0,1]$, but extends to $U_-\cup U_+$ by analytic
continuation. If $h$ denotes the lifting of $H$ to the universal
cover of the disk $D_*(0, R)$ by $\exp$, then we obtain {\em Abel's  
functional equation} 
\[ h(G(z)) = h(z) - \log \tau\] 
which allows one to interpret $h$ as the Fatou coordinate and  $U_{\pm}$
as the petals of $G$ at $x_0$. 
It also shows the nature of the
singularity of $H$ at $x_0$. Since the Fatou coordinate is 
$\log H = 
C_0 (z-x_0)^{-2} + C_1 (z-x_0)^{-1} + C_2 \log (z-x_0) + O(1)$, $C_0 < 0$, 
we get 
\[ H(z) = (z-x_0)^{C_2} \exp (\frac{C_0}{(z-x_0)^2}+\frac{C_1}{z-x_0}) 
\exp(\phi(z)) \]
where $\phi(z)$ is holomorphic.
The flat exponential factor
precludes $H$ from being analytic at $x_0$.

\section{Limits as $\ell_m\to \infty$}
In this section we prove Theorem~\ref{theo:10fa,1}.
\subsection{Bounds} 
\paragraph{Real bounds.}
For all results of this section, we assume that a unimodal mapping
$H_{\ell}(x) = |E_{\ell}(x)|^{\ell}$ is given, 
infinitely renormalizable with a periodic
combinatorial pattern $\aleph$, and satisfying the functional
equation~(\ref{equ:10fa,2}) with some scaling factor $\tau_{\ell}$.  

\begin{prop}\label{rb}
For every combinatorial pattern $\aleph$ there exist
two constants $1<T_1<T_2<\infty$,
such that $T_1<\tau_\ell<T_2$, for all $H_\ell$.
\end{prop}

The proof is contained in the following two lemmas~\ref{T_1},
~\ref{lem:27xp,1}.

First, let's make the following comment.
Given a solution $H_\ell(x)=|E_\ell(x)|^\ell$
of the equation~(\ref{equ:10fa,2})
with the constant $\tau=\tau_m$, let's
introduce a map $g(x)=E_\ell(|x|^\ell)$. Then
$g(0)=1$, $g$ is an even map, 
and $0$ is the critical point of the
unimodal map $g:[-1,1]\to [-1,1]$. It satisfies the
fixed-point equation 
\begin{equation}\label{equ:ggggg}
 \alpha g^{|\aleph|}(x) = g(\alpha x) \; ,
\end{equation}
where $\alpha=\alpha_\ell$ is a constant, which is either
$+\sqrt[\ell]{\tau}$ or $-\sqrt[\ell]{\tau}$.
Vice versa, to every solution $g(x)=E_\ell(|x|^\ell)$,
of the equation~(\ref{equ:ggggg}),
where $E_\ell$ is a diffeomorphism, there corresponds a solution
$H_\ell(x)=|E_\ell(x)|^\ell$ of ~(\ref{equ:10fa,2}) with
$\tau=|\alpha|^\ell$. 
One should have in mind the
following identity between $H$ and first return maps
of $g$ near the critical value $g(0)=1$ of $g$: 
\begin{lem}\label{lem:id}
For every $n\ge 0$,
\begin{equation}\label{equ:iden}
H_\ell(x)=\Lambda_n^{-1}\circ g^{|\aleph|^n}\circ \Lambda_n(x) \; ,
\end{equation}
where $\Lambda_n(x)=E_\ell(\tau^{-n} x)$ is a diffeomorphism
of $[0,1]$ onto its image.
\end{lem}
\begin{proof}
For $x\in [0,1]$,
one can write:
$\Lambda_n^{-1}\circ g^{|\aleph|^n}\circ \Lambda_n(x)=
\Lambda_n^{-1}\circ g\circ g^{|\aleph|^n-1}\circ \Lambda_n(x)=
\Lambda_n^{-1}\circ g\circ g^{|\aleph|^n-1}\circ g(|\alpha^{-n} x^{1/\ell}|)
=\Lambda_n^{-1}\circ g\circ g^{|\aleph|^n}(|\alpha^{-n} x^{1/\ell}|)
=\Lambda_n^{-1}\circ g(\alpha^{-n}g(|x|^{1/\ell}))=
\tau^n E_\ell^{-1}\circ E_\ell(|\alpha^{-n}g(|x|^{1/\ell})|^\ell)=
|g(|x|^{1/\ell})|^\ell=H_\ell(x)$.
\end{proof}

\begin{lem}\label{T_1}
There exists $1<T_1$,
such that $T_1<\tau_\ell$ for all $H_{\ell}$.
\end{lem}
\begin{proof} This follows easily from real bounds of~\cite{LvS}.
Indeed, let $U_n$ be the central $p^n$-periodic
interval of $g_\ell$, so that the endpoints of
$U_n$ are $u_n, -u_n$, where $u_n$ is $p^n$-periodic
point
of $g_\ell$. By the functional equation,
$u_n=u_0/\alpha_{\ell}^n$, where $u_0<-1$ is a fixed point
of $g_\ell$.
Let $I\supset g_\ell(U_n)$ be the maximal interval
on which $g_\ell^{p^n-1}$ is monotone. Then
$g_\ell^{p^n-1}(I)$ is contained in $U_{n-1}$.
On the other hand, by~\cite{LvS} (Lemma 9.1+Sect. 11),
there exists a universal constant $C_0$, such that
each component of $g_\ell^{p^n-1}(I)\setminus U_n$
has length at least $C_0|U_n|/\ell$ provided $n$ is large enough. 
Therefore,
$|u_0\alpha_{\ell}^{-n+1}/(u_0\alpha_{\ell}^{-n})|\ge 1+C_0/\ell$,
i.e. $|\alpha_{\ell}|> 1+C_0/\ell$,
and the existence of the universal $T_1$ follows.
(Let us remark that all real bounds of~\cite{LvS}
and their proofs hold without any changes for every
unimodal map of the form $E(|x|^\ell)$ where 
$E$ is a diffeomorphism of the Epstein class and
$\ell>1$ is any real number.)
\end{proof}

\begin{lem}\label{lem:27xp,1}
For every combinatorial type $\aleph$ there exists 
$T_2$ such that for all $H_{\ell}$ with combinatorial type $\aleph$, 
we get $\tau_{\ell} < T_2$. 
\end{lem}
\begin{proof} 
Decompose $H_{\ell} = |E_{\ell}|^{\ell}$. 
Let $(Z_1,Z_2)$ denote the maximal domain of monotonicity of
$H_{\ell}$ containing $0$ and $1$.  From~\cite{LvS}, 
$|E_{\ell}(Z_1)| \geq \sqrt[\ell]{\sigma}$ where $\sigma>1$ is
independent of $\ell$, though it might depend of $\aleph$.

A key
estimate here follows Lemma 3.8 in~\cite{bruke} and can be stated as
follows. Choose $0<A<1$ 
and let $B=H_{\ell}(A)$ (which is necessarily positive). 
Let us estimate from above the
$|H'_{\ell}(A)|$.  Consider the infinitesimal cross-ratio formed
by points $T, 0, A, A+dx$ where $Z_1 \leq T < 0$ is chosen so that 
$E_{\ell}(T) = \sqrt[\ell]{\sigma}$.   
Since $E_{\ell}$ is in the Epstein class, the
cross-ratio inequality gives
\[ |E_{\ell}'(A)| \frac{|t - 1|}{T}
\frac{A}{|b-1|} \frac{|A-T|}{|b-t|}  < 1\; \]
where we denoted $b=\sqrt[\ell](B)$ and $t=E_{\ell}(T)$.  
Since $\frac{|A-T|}{|T|} > 1$ and $|H'_{\ell}(A)| = \ell
b^{\ell-1} |E_{\ell}'(A)|$, we get
\[ |H'_{\ell}(A)| < \frac{|b-t|}{|t-1|} \frac{|b-1|}{A}
\frac{B}{b} \ell = \frac{B}{A} \ell |b-1| \ell \frac{|t-b|}{t}
\frac{1}{\ell |t-1|} \frac{t}{b} \; .\]
Since $|b-1| < \log b^{-1}$, $\frac{|t-b|}{t} < \log \frac{t}{b}$ 
$|t-1| > \log t$, we get 
\[ |H'_{\ell}(A)| < \frac{B}{A} \log b^{\ell} \log
\frac{t^{\ell}}{B} \frac{1}{\log t^{\ell}} \frac{t}{b} \; .\]
Finally recalling that $t= \sqrt[\ell]{\sigma}$ and
$b^{\ell} = B$, 
we get 

\begin{equation}\label{equ:9xp,1}
 |H'_{\ell}(A)| < \frac{B}{A} \log B^{-1} \log
\frac{\sigma}{B} \frac{1}{\log \sigma} 
\sqrt[\ell]{\frac{\sigma}{B}}\; .
\end{equation}

When $A=1, H_{\ell}(1),\cdots, H_{\ell}^{|\aleph|-2}(1)$ then 
$B=H_{\ell}(A)$ is at least $\tau^{-1}_{\ell}$, since
$\tau_{\ell}^{-1} = H_{\ell}^{|\aleph|-1}(1)$ is the closest return of
the orbit of $0$ to itself. For all such $A$, we can thus
rewrite~(\ref{equ:9xp,1}) as 

\[  |H'_{\ell}(A)| < \frac{B}{A} \log \tau{\ell}\frac{\log \tau_{\ell} +
1}{\log\sigma} \sqrt[\ell]{\frac{\sigma}{B}}\; . \]

Now the functional equation implies that 
$|(H_{\ell}^{|\aleph|-1})'(1)| = 1$. 
Therefore, if we take  the product of such estimates for all 
$A$ equal to 
\[ 1, H_{\ell}(1),\cdots, H_{\ell}^{|\aleph|-2}(1) \; ,\]
we get $1$ on the left-hand side. 

We obtain 
\[ 1 < \tau_{\ell}^{-1} 
\left(
\log\tau_{\ell}\frac{\log\tau_{\ell}+1}{\log\sigma}
\right)^{|\aleph|-1} \times \]
\[ \times\left(\tau_{\ell}\sigma \right)^{(|\aleph|-1)/\ell} \;.\]
Since for $\ell > |\aleph|$ the right-hand side goes to $0$ as
$\tau_{\ell}$ increases to $\infty$, the
estimate follows for all $\ell$ but finitely many.  
\end{proof}

\paragraph{Complex bounds.}

\begin{prop}\label{cb}
For every combinatorial type $\aleph$, there exist constants
 $\ell_0$, $\lambda>1$ and $R_1$, such that, 
for every $H_{\ell}=|E_\ell|^\ell$ with combinatorial type $\aleph$ 
which satisfies the functional
equation~(\ref{equ:10fa,2}) with some $\tau_{\ell} > 1$, as soon as  
$\ell\geq \ell_0$, there exists $1
< R < R_1$ as follows.
The function $E_\ell$ extends to a map
from the Epstein class defined on a neighborhood of $[0,1]$, 
so that 
function $H_{\ell}=|E_\ell|^\ell$ extends to a unimodal
function from some interval $[R'_-, R'_+]$ onto
$[0,R]$, having a fixed point $b_{\ell} \in (1, R)$, 
with the following inequalities:
\[ |R'_-| \leq |R'_+| \leq \lambda^{-1} R\; .\]
\end{prop}

The name ``complex
bound'' comes from the fact that since for $\ell$ which is an even
integer $H_{\ell} = (E_{\ell})^{\ell}$ with $E_{\ell}$ in the Epstein
class, Proposition~\ref{cb} implies that $H_{\ell}$ has a
polynomial-like extension onto the domain $D(0, R)$.
\begin{proof}
Proposition~\ref{cb} follows from~\cite{LvS}. 
To make the reduction, we consider the dynamics of the corresponding
map $g(x)=E_\ell(|x|^\ell)$ on the level
of $p^n$-periodic central interval $U_n$
where $n$ is large enough. To connect this dynamics
with the map $H$, one can use the identity(~\ref{equ:iden})
rewritten in the form
$E_\ell^{-1}(x)=\tau^{-n}E_\ell^{-1}\circ g^{-(p^n-1)}(|\alpha|^{-n}x)$,
where $g^{-(p^n-1)}$ is the branch from the interval $U_n$
to a neighborhood of $g(0)=1$. The identity holds
originally in a small neighborhood of $0$. On the other hand,
the right-hand side extends to a real-analytic function 
on an interval $[-R^{1/\ell}, R^{1/\ell}]$, where 
$R=|\alpha^n g^{p^n}(\tilde u)|^\ell$ and $\tilde u$
is a point defined in Lemma 9.1 of~\cite{LvS}
for the $p^n$-periodic central interval
of $g$. Then we apply the latter Lemma and get the result.
\end{proof}

\subsection{Limit maps}
Our aim is to pick a convergent subsequence from 
$H_{\ell_m}$ by some kind of compactness argument. 
The problem is that as $\ell_m\to \infty$,
then the domains of definition as $U_{\ell_m}$ 
tend to degenerate at a limit of the critical points $x_{\ell_m}$.

To deal with this phenomenon, we consider inverse branches of $H_{\ell_m}$
corresponding to values to the left and to the right of
the point $x_{\ell_m}$.

From the form and normalization of mappings $H_{\ell_m}$, each
of them can be represented as $|E_{\ell_m}(x)|^{\ell_m}$ with
$E_{\ell_m}$ an Epstein diffeomorphism mapping at least onto the
interval $(-\sqrt[\ell_m]{R_{\ell_m}}, \sqrt[\ell_m]{R_{\ell_m}})$
with $R_{\ell_m}$ chosen from Proposition~\ref{cb}. Further from
Proposition~\ref{cb} one gets that $E_{\ell_m}^{-1}(D(0,
\sqrt[\ell_m]{R_{\ell_m}})) \subset D(0, \lambda^{-1} R_{\ell})$. 
By taking a subsequence we can assume without loss of generality that 
$R_{\ell_m} \rightarrow R\geq \lambda > 1$. Similarly, in the light of
Proposition~\ref{rb}, we may assume that $\tau_{m} \rightarrow
\tau > 1$. Choosing yet another subsequence, we may assume that
$x_{\ell_m} \rightarrow x_0$.   

We will actually invert not $H_{\ell_m}$, but its lifting $h_{\ell_m}$
to the universal cover of $D_*(0, R)$ by $\exp$. 
This will have two real branches, one mapping onto a
right neighborhood of $x_{\ell_m}$ and one onto a left neighborhood. 
Their complex extensions are 
\begin{equation}\label{equ:12ga,1}
P_{\ell_m}^+(w) :=
E_{\ell_m}^{-1}(\exp(w/\ell_m))
\end{equation}
\[ P_{\ell_m}^-(w) := E_{\ell_m}^{-1}  (-\exp(w/\ell_m))\; .\]

 Both transformations 
are defined in $\Pi_{m}
:= \{ w :\: \Re w < \log R_{\ell_m} \}$ and map into $D(0, \lambda
R_{\ell_m})$ by Proposition~\ref{cb}. 

By Montel's theorem we can  pick a subsequence $m_k$, such that 
$P^{\pm}_{\ell_{m_k}}$ converge to mappings $P^{\pm}$
defined on $\Pi_* := \{ w :\: \Re w < \log R\}$. Since the domains
vary with $m$, they should be normalized for example by precomposing
with a translation, which tends to $0$ in the limit. This implies 
uniform convergence on compact subsets, with the understanding that every
compact subset of $\Pi_*$ belongs to $\Pi_{m}$ for almost all $m$.        
In the sequel, we will ignore this subsequence and simply assume that 
$P^{\pm}_{\ell_m}$ converge.

Let us see that $P^{\pm}$ are both non-constant.
Note that $P^{+}_{\ell_m}(0)=0$.
Moreover, by the functional equation,
$H_{\ell_m}^{p+1}(x_{\ell_m})=H_{\ell_m}^p(0)=1/\tau_{m}$, and,
by the combinatorics,
$H_{\ell_m}(1/\tau_{m})=H_{\ell_m}^p(1)\in (1/\tau_{m}, 1)$.
Therefore,
there exists a point $a_{\ell_m}\in (\log(1/T_2), 0)$, such that
$P^{+}_{\ell_m}(a_{\ell_m})=1/\tau_{m}\subset (1/T_2, 1/T_1)$,
so that $P^{+}_{\ell_m}(a_{\ell_m})$ are uniformly away from
zero. Similarly, one can see that 
any limit function of the family $\{P^{-}_{\ell_m}\}$
is not constant as well.
The considerations are slightly different in the cases
$p=2$ and $p>2$; for example, let $p>2$. Then
$H_{\ell_m}^2(0), H_{\ell_m}^{2+p}(0)\in (H_{\ell_m}^p(0),1)=
(1/\tau_{m},1)\subset (1/T_2,1)$; on the other hand, 
$P_{{\ell_m}}^+(a_{{\ell_m}})=1/\tau_{m}$,
where $a_{\ell_m}=\log(H_{\ell_m}^{p}(1))\in (\log(x_{\ell_m}),0)$;
the limit maps of $(P_{\ell_m}^+)$ are not constants,
hence, there is $c^*<1$ such that $H_{\ell_m}^p(1)<c^*$ for all ${\ell_m}$;
therefore, $P_{\ell_m}^-(\log(H_{\ell_m}^{1+p}(1)))=H_{\ell_m}^p(1)<c^*$
while $P_{\ell_m}^-(\log(H_{\ell_m}(1))=1$, and the conclusion follows. 

It is also clear that $P^{\pm}$ are both univalent. This is because
for any compact subset of $\Pi_*$ and $\ell_m$ large enough,
$P^{\pm}_{\ell_m}$ are univalent on this set, which is evident from their
defining  formulas~(\ref{equ:12ga,1}).

Let us define $x_0^{\pm} : = \lim_{x \rightarrow -\infty} P^{\pm}(x)$.  
Since $(P^+)^{-1}$ in increasing on $(x_0^+, 1]$ and $(P^-)^{-1}$ is
decreasing on $[0,x_0^-)$, we must have $x_0^- \leq x_0^+$. We will
next show that $P^+(\Pi_*) \subset {\cal D}((x_0^+, R'_+),\pi/2)$ and 
$P^-(\Pi_*) \subset {\cal D}((x_0^-, R'_-),\pi/2)$. We used here
notations $R'_\pm$ from the statement of Proposition~\ref{cb} and for
any interval $I$, ${\cal D}(I, \pi/2)$ means the Euclidean disk with
$I$ as its diameter. We will concentrate on the first inclusion. It
will follow once we show that for $m$ large enough and any $w\in \Pi$,
$\exp(w/\ell_m) \in {\cal D}(0, \sqrt[\ell_m]{R})$, by
formula~(\ref{equ:12ga,1}) and since $E_{\ell_m}$ is in Epstein class.  
The inclusion follows since 
$ |\arg(\log R - w)| < \pi/2$ and $\exp$ is
conformal, so 
\[ \lim_{m\rightarrow \infty} \arg (\sqrt[\ell_m]{R} - \exp(w/\ell_m))
= \arg(\log R - w) \; .\]

\paragraph{Checking conditions for the Eckmann-Wittwer class.}
We can now define a limit mapping $H$ which will be shown to satisfy
Definition~\ref{defi:10fa,1}. 

We set $U_\pm = P^{\pm} (\Pi_*)$. Then $H_{|U_\pm} :=
\exp\circ (P^{\pm})^{-1}$. $H$ can also be defined and equal to $0$ 
on the interval (perhaps degenerate) $[x_0^-, x_0^+]$. 

We have shown that
$H_{\ell_m}$ converge to $H$ uniformly on any compact subset of 
$(x_0^+,R'_+]$ or $[R'_-,x_0^-)$, again using notations from
Proposition~\ref{cb}.  Because the mappings $H_{\ell_m}$
and $H$ are unimodal, this implies uniform convergence on compact
subsets of $(R'_-, R'_+)$. 

Setting out to check the conditions of Definition~\ref{defi:10fa,1},
we see that the functional equation is satisfied simply by passing to
the limit with $m$. In particular, we use the fact that since
$H_{\ell_m}$ converge uniformly, their family is
equicontinuous. 

The conditions second, third and fourth are satisfied by construction. 

To derive the fifth condition, observe that $H(1) < 1$ while $H(R'_+)
= R > R'_+$. So, there must be a fixed point $b_0$ between $1$ and
$R'_+$ which is unique and repelling because $H$ has non-positive
Schwarzian derivative in the light of condition 4.

With regard to the sixth condition, the symmetry with respect to the
real line follows from
formulas~(\ref{equ:12ga,1}). We have proved the disjointness of the
closures of $U_-$ and $U_+$ except perhaps if $x_0^+ = x_0^-$. So we
now need to prove this equality. This will require another idea and we
will in fact prove property 7 first.

\paragraph{The associated dynamics of $G$.}
For every $m$, define $G_{\ell_m}(z)=H_{\ell_m}^{p-1}(z/\tau_{m})$
which is well-defined and holomorphic
in a neighborhood of the point $x_{\ell_m}$. The functional equation
yields 
\begin{equation}\label{equ:13gp,1}
\tau^{-1}_{m} H_{\ell_m} = H_{\ell_m} \circ G_{\ell_m}
\end{equation}
 on the interval $[0,1]$. Since $H_{\ell_m}(x) = 0$ implies $x =
x_{\ell_m}$, the functional equation implies that $x_{\ell_m}$ is a
fixed point of $G_{\ell_m}$. Since $|H_{\ell_m}(x_0+x)| = A
|x|^{\ell_m} + o(|x|^{\ell_m})$, expanding $G$ into the power series
and substituting into~(\ref{equ:13gp,1}) yields $|G'((x_0)| =
\tau_{m}^{-1/\ell_m}$. Also, equation~(\ref{equ:13gp,1}) and the
fact that $H_{\ell_m}$ is unimodal imply that $x_{\ell_m}$ attracts the
entire interval $[0,1]$ under the iteration of $G_{\ell_m}$.

Since the fixed point equation remains valid for the limit function
$H$, if we define $G(x) = H^{p-1}(\tau^{-1}x)$,
equation~(\ref{equ:13gp,1}) is also satisfied with indices $\ell_m$
removed. We see that $G(x_0) = x_0$ and $x_0$ is topologically
non-repelling: $|G(x) - x_0| \leq |x-x_0|$ for every $x\in [0,1]$.  

Recall now that $H^{-1}(0) = [ x_0^+, x_0^- ] \ni x_0$.  
\begin{lem}\label{lem:13gp,1}
$G([x_0^-, x_0^+]) = [x_0^-, x_0^+]$. 
\end{lem}
\begin{proof}
From the functional equation, since $\tau^{-1} H([x_0^-,x_0^+]) = 0$,
it follows that $G([x_0^-, x_0^+]) \subset [x_0^-, x_0^+]$. If it were
a proper subset however, we would have $G(x) \in [x_0^-, x_0^+]$ for
some $x\notin [x_0^-, x_0^+]$, which would imply $H(x) = 0$ contrary
to $[x_0^-, x_0^+] = H^{-1}(0)$. 
\end{proof}

\begin{lem}\label{lem:13gp,2}
On a neighborhood of the interval $[0,1]$ in the complex
plane $G(z) = H^{p-1}(z/\tau)$ is well defined, in particular analytic.
\end{lem}
\begin{proof}
Denote $K=[0, \tau^{-1}]$.
To show the claim of the lemma, it is enough to show 
that $H^n(K) \cap [x_0^-, x_0^+]=\emptyset$ 
for any
$0\leq n\leq p-2$. 
Otherwise, for some $0\leq j\leq p-2$, $0\in H^{j+1}(K)$.
On the other hand, $K=[H(x_0), H^{p+1}(x_0)]$, hence,
by the combinatorics, the intervals $H^n(K)$, $0\leq n\leq p-2$,
are pairwise disjoint, a contradiction.
\end{proof}

From Lemma~\ref{lem:13gp,2} we conclude
that $G_{\ell_m}$ converge to $G$ uniformly on a complex
neighborhood of $[0,1]$ and that $G$ restricted to $[0,1]$ is a
diffeomorphism in the Epstein class, in particular $SG \leq 0$. Since
$|G'_{\ell_m}(x_{\ell_m})| = \sqrt[\ell_m]{\tau^{-1}_{m}}$, the
convergence implies $(G^2)'(x_0) = 1$. Coupled with the information
that $x_0$ is topologically non-repelling on both sides, this implies
the power-series expansion: 
\[ G^2(z) - x_0=(z-x_0)+a(z-x_0)^{q+1}+O(|z-x_0|^{q+1})\]
with some $a\leq 0$ and some $q$ even. 
First, we prove that $a\not=0$, i.e. $G^2$ is not the identity.
If $G^2(z)=z$, then, for every $x\in [0,1]$,
$H(x)=H(G^2(x))=H(x)/\tau^2$, i.e. $H(x)=0$ and
$[x_0^-,x_0^+]=[0,1]$, a contradiction. Thus, $a<0$.

Now we prove that $q=2$ considering a perturbation.
There is a fixed complex neighborhood
$W$ of $x_0$,
such that the sequence of maps $(G^2_{\ell_m})^{-1}$
are well-defined in $W$ and converges 
uniformly in $W$ to $(G^2)^{-1}$.
Since each $H_{\ell_m}$ belongs to the Epstein class,
then each $(G^2_{\ell_m})^{-1}$
extends to a univalent map
of the upper (and lower) half-plane
into itself. It extends also continuously on the real line,
and has there exactly one fixed point, which is
$x_{\ell_m}$ and which is repelling. Therefore, by the Wolff-Denjoy theorem,
$(G^2)^{-1}$ has at most one fixed point in either half-plane, and one
which is strictly attracting. Thus, for any $m$, $G^2_{\ell_m}$ has at
most three simple fixed points on $W$, which implies $q=2$ by Rouche's
principle.  In this way, we have proved condition 7. 

Finally, we can finish the proof of condition 6 by showing that
$x_0^-=x_0^+=x_0$. Indeed, $x_0^-$ and $x_0^+$ are both fixed points of
$G^2$ by Lemma~\ref{lem:13gp,1}, but the local form $G$ at $x_0$ and
the condition $SG\leq 0$ mean that $x_0$ is the unique fixed point of
$G$ on $[x_0^-, x_0^+]$. 

We have finished the proof of Theorem~\ref{theo:10fa,1}.

\section{Dynamics of EW-maps}\label{sec:15fa,1}
In this section, we will construct basic dynamical theory of EW-maps,
including the construction of their Julia sets and quasiconformal
equivalence.  
\subsection{Real dynamics}
Recall that an interval is called {\em wandering} for a unimodal map 
provided that all its forward images avoid the critical point and its
$\omega$-limit set is not a periodic orbit. 

\begin{prop}\label{prop:10fp,1}
If $H$ is a mapping in the EW-class with any combinatorial pattern
$\aleph$, then $H$ has no wandering interval.  
\end{prop}

Set $p := |\aleph|$ and let $I_0 =(b'_0, b_0)$ using the notation of 
Definition~\ref{defi:10fa,1}. 
 
We have  the functional identity $H^{p^n}=G^n\circ H\circ G^{-n}$ for
any $n$ on $I_n := G^n(I_0)$. To verify the identity, act on both sides by
$G^n$ from the left and use the functional equation $H\circ G =
\tau^{-1} H$ and the definition $G(x) = H^{p-1}(x\tau^{-1})$ $p$
times.

Then $G^m$ provides a smooth conjugacy
between $H^{p^m}$ on $I_m$ and $H$ on $I_0$. Since for either
connected component $C$ of $I_0 \setminus I_1$ intervals $C,\cdots,
G^{m-1}(C)$ belong to $I_0$ and 
are pairwise disjoint, the distortion of $G^m$ on $C$ is
bounded in terms of the total nonlinearity of $G$ on $I_0$ and
independently of $m$.

Introduce the following sets of intervals:
for every $m\ge 1$, let $\{I_{m,j}\}_{j}$ be the collection
of all connected components of the first entry map from $I_{m-1}$ into
$I_m$. These intervals cover $I_{m-1}$ except for countably many
points (preimages of the endpoints of $I_m$). 
Define dynamics $F$ on $P=\cup_{m\ge 1}\cup_{j} I_{m,j}$:
if $x\in P_m=\cup_{j} I_{m,j}$, 
then $F(x)=H^{p^{m-1}}(x)$.
Then
$F$ maps homeomorphically any $I_{m,j}$ onto another interval 
$I_{m,j'}$ and eventually onto $I_m$.

Let $\rho_A$ denote the hyperbolic metric
on an interval $A=(a,b)$, i.e. 
\[ \rho_A(x,y) = | \log\frac{|x-a| |b-y|}{|y-a| |b-x|} |\; \]
and  denote by $\rho_P$ the metric
on $P$, defined so that $\rho_P(x,y) = \rho_{I_{m,j}}$ if $x,y \in
I_{m,j}$ or is $\infty$ if no such $m,j$ exist.

Start with following lemma.

\begin{lem}\label{lem:10fp,1}
There exists a constant $K$ such that for every $m, j$ the length of 
$I_{m,j}$ in $\rho_{I_{m-1}}$ is less than $K$. 
\end{lem}
\begin{proof}
Because $G^{-m}$ maps intervals $I_{m,j}$ onto $I_{0,j}$ and every
$I_{0,j}$ is contained in a connected component of $I_0 \setminus I_1$,
and because of uniformly bounded distortion, without loss of
generality we can set $m=0$. If a sequence $j_k$ exists such that the
lengths of $I_{0,j_k}$ go to $\infty$ then, perhaps by taking a
subsequence, right endpoints of $I_{0,j_k}$ tend to $b_0$. But if
$I_{0,j_k} = (\alpha,\beta)$ with $\beta$ close to the repelling fixed
point $b_0$, then $\alpha > H(\beta)$ since $(H(\beta), \beta)$
contains a preimage of an endpoint  of $I_1$. Thus, the hyperbolic
length of $(\alpha,\beta)$ can be bounded in terms of the eigenvalue
of $H$ at $b_0$. 
\end{proof}

Supposing now that a wandering interval $J$ exists, we observe that 
$J$ must be disjoint from $\partial I_m$ for every $m$. This is
because the endpoints of $I_m$ are pre-repelling fixed points of
$H^{p^{m-1}}$ and every one-sided neighborhood of such a fixed point
will eventually cover $x_0$ under the iteration of
$H^{p^{m-1}}$. Then, $J \subset I_{n_0, j_0}$ for some $n_0, j_0$. 
Consider the sequence of intervals $J_k:=F^k(J)$
Then $(\rho_P(J_k))_{k\ge 0}$ is an increasing sequence.
Moreover, each time $J_k$ is mapped into $I_m$ for the first time, 
$\rho_P(J_k) > \lambda \rho_P(J_{k-1})$. Here $\lambda$ is the
expansion constant of the inclusion map $I_{m,j} \rightarrow I_m$,
where $J_k \subset I_{m,j}$, with the metric $\rho_{P} =
\rho_{I_{m,j}}$ in the domain and $\rho_{I_m}$ in the image. 
Observe that  $\lambda$ is bounded
away from $1$ by Lemma~\ref{lem:10fp,1}.

Hence, $\rho_P(J_k)$ goes to $\infty$ with $k$. But as soon as $J_k
\subset I_m$, then $J_k$ is also wandering for $H^{p^m}$ and so
contained in $I_{m,j}$, which leads to a contradiction with 
Lemma~\ref{lem:10fp,1}.

\subsection{Julia set}
Recall that for any
interval $I$ and $0<\theta<\pi$ the set ${\cal D}(I,\theta)$ consists
of all points in $\CC$ whose distance to the ``line'' $I$ in the
hyperbolic metric of $(\CC \setminus \RR) \cup I$ is less than a
constant. Such a set is bounded by arcs of circles which intersect 
$\RR$ at the endpoints of $I$ and $\theta$ denotes the angle formed by
these arcs with $\RR$ with the convention that ${\cal D}(I,\theta)$
grows with the growth of $\theta$, see~\cite{destrien} and~\cite{Su}.

\paragraph{A few lemmas.}
We begin with couple of lemmas describing the complex dynamics of 
$H$. 

\begin{lem}\label{lem:28fp,1}
Let $H$ belong to the EW-class with some combinatorial type $\aleph$. For
every $n = 0,\cdots$ consider real points $u_-^n < x_0 < u_+^n$
defined by $H(u_-^n) = H(u_+^n) = \tau^{-n} R$. Consider a point $z\in \CC$
and $k=1,\cdots$ chosen so that 
\[ H^k(z) \in {\cal D}((u_-^n, u_+^n), \pi/2)\] 
but $H^{k'}(z) \notin D(0,\tau^{-n} R)$ for all $0<k'\leq k$. 

For any such choice of $z,k,n$ there is an inverse branch of $H^k$
defined on ${\cal
D}((u_-^n,u_+^n), \pi/2)$ which sends $H^k(z)$ to $z$.  
\end{lem}
\begin{proof}
Since the Poincar\'{e} neighborhood is simply connected, the only
obstacle to constructing the inverse branch may be if the omitted
value $0$ is encountered. Thus suppose that for some $k' > 0$,
$\zeta$, which is an inverse branch of $H^{k-k'}$ well defined on 
${\cal D}((u_-^n, u_+^n), \pi/2)$, maps
$H^k(z)$ to $H^{k'}(z)$ and its image contains $0$. It follows that 
$H^{k-k'}(0) \in (u_-^n, u_+^n)$ and so $H^{k-k'+1}(0) < \tau^{-n} R$. It
follows that $k-k'+1$ must be a multiple of $p^n$, where we denote
$p:=|\aleph|$. Then $\zeta$ is just a real map on the real line and 
$\zeta(u_-^n, u_+^n) \subset [0,\tau^{-n} R)$. But since Poincar\'{e}
neighborhoods are mapped into Poincar\'{e} neighborhoods of the same
angle by $\zeta$, we get $H^{k'}(x) \in \zeta({\cal D}(u_-^n,u_+^n))
\subset \D(0,\tau^{-n} R)$ contrary to the hypothesis of the lemma. 
\end{proof}

\begin{lem}\label{lem:28fp,2}
Let $H$ belong to the EW-class. Define $U_{+,c}$ to be the connected
component of $H^{-1}(D(0,R) \cap \{z: \Re(z)>0\})$ which contains $U_+\cap
\RR$. Also, specify $U_{-,c}$ analogously. For some point $z\in \CC$ 
suppose that $H^k(z)\in U_{+,c}\cup U_{-,c}$ for all $k\in \NN$ and
the Euclidean distance from the forward orbit of $z$ to the
$\omega$-limit set of $0$ is $0$. Then $z \in \RR$. 
\end{lem}
\begin{proof}
First, we observe that $0$ must belong to the closure of the forward
orbit of $z$. Indeed, by hypothesis, the orbit of $z$ is contained in
$H_{+,c} \cup H_{-c} \cup \{x_0\}$ and $H$ restricted to this set is
continuous. Then by the minimality of the $\omega$-critical set, if
the orbit of $z$ accumulates on it somewhere, then it also accumulates
at $0$. As soon as $0$ belongs to the $\omega$-limit set of $z$, we can find a
sequence of iterates $k_n$, perhaps not strictly increasing, such that
$H^{k_n+1}$ are first entry times of $z$ into $D(0, \tau^{-n} R)$.
Then $H^{k_n}(z)$ belong to ${\cal D}((u^n_-, u^n_+), \pi/2)$ by the
Epstein class properties postulated in Definition~\ref{defi:10fa,1}.  
Consequently,  we can consider inverse branches $\zeta_n$ constructed
by Lemma~\ref{lem:28fp,1}. Since the orbit of $z$ is contained in
$U_{+,c} \cup U_{-,c}$, then each $\zeta_n$ will map $(u_-^n, u_+^n)$
into some real interval $T_n$ and $z \in {\cal D}(T_n, \pi/2)$. But
the lengths of $T_n$ have to go to $0$ or we could find a non-trivial
interval contained in infinitely many of them. Such an interval would
be wandering in contradiction to Proposition~\ref{prop:10fp,1}. It
follows that the distance from $z$ to $\RR$ must be $0$. 
\end{proof}

\paragraph{The filled-in Julia set.}

\begin{defi}\label{defi:29fa,1}
If $H$ belong to the EW-class we define its {\em filled-in Julia set} 
$K_{H}$ as follows:
\[ K_{H} := \{ z :\: \forall n\geq 0\: H^n(z) \in \overline{U^+ \cup
U^-\}}\ \cup_{n\ge 0}H^{-n}(\{x_0\}) ; .\]
\end{defi}

The disadvantage of Definition~\ref{defi:29fa,1} is that $K_H$ appears
to depend on the parameter $R$ from
Definition~\ref{defi:10fa,1}. Also, other than the name there is a
priori no connection between $K_H$ and Julia sets of globally defined
holomorphic mappings, so any theory has to be developed from scratch.  

\begin{theo}\label{pre}
For an EW-map $H$, the filled-in Julia set $K_H$ is the closure of the
set of all preimages of $0$ by iterates of $H$. In particular, $K_H$ is
independent of the particular choice of $R$ in
Definition~\ref{defi:10fa,1} and its interior is empty. 
\end{theo}

In the course of the proof we introduce some ideas
which will be used also later on.
Start by observing that $K_H \cap \RR = [b'_0, b_0]$ because of the
negative Schwarzian of $H$. On the other hand, preimages of $0$ are
dense in $[b'_0, b_0]$ in the light of
Proposition~\ref{prop:10fp,1}. Also, $x_0$ is not an interior point
of $K_H$ since it lies on the boundary of the domain of definition, so
once we know that the preimages of $0$, hence of $x_0$, are dense in
$K_H$, then $K_H$ indeed has a vacuous interior. So we only need to
prove the density of the set of preimages in $K_H$. This is done by
considering the hyperbolic metric. 

\paragraph{Hyperbolic metric.}
Let $\omega$ denote the $\omega$-limit set
of the critical point $x_0$ by the map $H:[0,1]\to [0,1]$;
The set $\omega$  is closed and forward invariant; moreover,
the set $V\setminus \omega$ is open and connected.  
Denote by $\rho$ the hyperbolic metric
of the domain $V\setminus \omega$.

If $\rho$ is a metric and $F$ a function, we will write $D_{\rho}
F(z)$ for the expansion ratio with respect to the metric $\rho$, thus 
\[ D_{\rho} F(z) = |F'(z)| \frac{d\rho(F(z))}{d\rho(z)} \; .\]

By Schwarz's lemma, we have $D_{\rho} H(z) > 1$ for every $z \in U_+
\cup U_-\setminus H^{-1}(\omega)$. We will prove that if $z\in K_H$ and no forward image of
$z$ is real, then 
\[ \lim_{n\rightarrow\infty} D_{\rho} H^n(z) = \infty\; .\] 

We will observe expansion of the hyperbolic metric based on the
following fact:  

\begin{fact}\label{fa:20gp,2}
Let $X$ and $Y$ be
hyperbolic regions and $Y\subset X$ and $z\in Y$. Let $\rho_X$ and
$\rho_Y$ be the hyperbolic metrics of $X$ and $Y$, respectively. 
Suppose that the hyperbolic
distance in $X$ from $z$ to $X\setminus Y$ is no more than $D$. 
For every $D$ there
is $\lambda_0 > 1$ so that  
$|\iota'(z)|_H\leq \frac{1}{\lambda_0}$, where $\iota :\: Y \rightarrow X$
is the inclusion, and the derivative is taken with respect to
the hyperbolic metrics in $Y$ and $X$, respectively.
\end{fact}

In our case, we will set $Y := V \setminus (\omega \cup H^{-1}\omega)$
and $X= V \setminus \omega$. It follows that $D_{\rho}H(z) \geq \lambda_d
> 1$ provided that the distance from $z$ to $H^{-1}(\omega)$ with
respect to $\rho$ is bounded by $d$. 

Fixing $z\in K_H$ which is not eventually mapped into $\RR$ and based
on Lemma~\ref{lem:28fp,2}, we distinguish two eventualities. The
first is that the Euclidean distance from the forward orbit of $z$ to
$\omega$ is positive. The hyperbolic
distance from $H^n(z)$ to $H^{-1}(\omega)$ is bounded uniformly in $n$
and $D_{\rho} H^n(z)$ grows at a uniform exponential rate. 

In the second (opposite) case, Lemma~\ref{lem:28fp,2} gives us a sequence $n_k$
such that $H^{n_k}(z) \notin U_{+,c} \cup U_{-,c}$. Now the hyperbolic
distance from $z' := H^{n_k}(z)$ to $H^{-1}(\omega)$ is uniformly
bounded. To see this, fix attention on the case when  $z' \in
U_+$. The hyperbolic metric $\rho_+$ on $U_+ \setminus \omega$ is
bigger than $\rho$ and $U_+$ can be conveniently uniformized by the
map $\log H$ where the branch of the $\log$ is chosen so that the map
is symmetric about the real axis. The image of $U_+$ is the half-plane 
$\{ w :\:  \Re w < \log R\}$, but $\Re \log H(z') < R'$ with fixed
$R'<R$ since otherwise $H(z') \notin \overline{U_+ \cup U_-}$.     
Since $z' \notin U_{+,c}$, then $|\Im \log H(z')| \geq \pi/2$. The set
$\log H (H^{-1}(\omega))$ is doubly periodic with periods $2\pi i$ and
$\log \tau$, so evidently the hyperbolic distance from $\log H(z')$ to
it is bounded.

Now suppose that $z\in K_H$ and $z$ is not in the closure of the set
of preimages of $0$. This implies that no forward image of $z$ is
real, so $D_{\rho} H^n(z)\rightarrow \infty$ as just argued. Moreover,
we have shown that for some sequence of iterates $H^{n_k}$, the
hyperbolic distance from $H^{n_k}(z)$ to $H^{-1}(\omega)$ is uniformly
bounded. By pulling back to $z$, we see that the hyperbolic distance
from $z$ to $\bigcup_{j=1}^{\infty} H^{-j}(\omega)$ is zero, and since 
$\omega$ is contained in the closure of the preimages of $0$, this
concludes the proof of Theorem~\ref{pre}.

\subsection{Quasi-conformal equivalence}
Let now $H:U^+\cup U^-\to V$ and $\hat H:\hat U^+\cup \hat U^-\to \hat V$
be two maps from the EW-class with the same combinatorial type $\aleph$. 
We will eventually show that $H=\hat{H}$, but as the first step, we
prove they are quasi-conformally conjugate. 

\begin{prop}\label{prop:11fp,1}
For every pair of maps $H, \hat{H}$, both in the EW-class with the same
combinatorial type $\aleph$, there exists a quasi-conformal homeomorphism 
$\phi_0$ of the plane, symmetric w.r.t. the real axis,
and normalized so that $\phi_0(0)=0, \phi_0(1)=1$,
which  conjugates $H$ and $\hat H$, i.e. $\phi_0(U_-) = \hat{U}_-$,
$\phi_0(U_+) = \hat{U}_+$ and  
$\phi_0\circ H(z)=\hat H\circ \phi_0(z)$ for every $z\in U_+ \cup
U_-$. 
\end{prop}

The proof of Proposition~\ref{prop:11fp,2} will be obtained from the
following 

\begin{prop}\label{prop:11fp,2}
For every pair of maps $H, \hat{H}$, both in the EW-class with the same
combinatorial type $\aleph$, there is a 
mapping $\phi_1$ defined and continuous in $\overline{U_- \cup U_+}$, 
quasi-conformal in
the interior, symmetric about the real axis, 
and which can be restricted to a quasi-symmetric
orientation-preserving map of the interval $\RR\cap(U_-\cup U_+)$.
Dynamically, $\phi_1(H(z)) = \hat{H}\phi_1(z)$ for every $z$ in the
forward orbit of $x_0$ and $\frac{\hat{R}}{R} H(z) =
\hat{H}(\phi_1(z))$ for every $z \in \partial (U_- \cup U_+)$.  
\end{prop}

Given Proposition~\ref{prop:11fp,1}, one can apply the pull-back
argument~\cite{Su}  to the original maps $H, \hat H$. Once
Theorem~\ref{pre} has been 
established, the construction becomes standard.

\paragraph{Presentation functions.}
In order to show Proposition~\ref{prop:11fp,2}, we have to find an
alternative to the standard method which consists of constructing
first a quasi-symmetric equivalence on the real line based on the
bounded geometry of the Cantor attractor.  
For maps in the EW-class, however, 
the $\omega$-critical set $\omega$ has no bounded geometry,
because $\omega$ is invariant under the map $x\mapsto G(x)$.
Instead, given $H$,
we construct a complex 
box mapping $h=h_H$ with simpler dynamics (post-critically finite to
be precise)
so that $\omega$ is a subset of ``repeller'' of such map. This
generalizes the idea of ``presentation functions'',
see~\cite{misled},~\cite{CEp222}, 
which was to realize a non-hyperbolic attractor as
hyperbolic repeller.

Write $p := |\aleph|$. Recall the notation $I_1 = G(I_0)$ from the
proof of Proposition~\ref{prop:11fp,1}. 
Introduce intervals $J_1 = (0, b_0\tau^{-1})$, $J_p := I_1$
and $J_q$ which is the connected component of $H^{q-p}(J_p)$ which
contains $H^{q-1}(J_1)$ for $1<q<p$. Then $J_q$, $q=1,\cdots,p$ are
pairwise disjoint intervals which cover $\omega$ and are contained in 
$(0, R')$ for some $R' < b_0$. Also, $H(J_p) = J_1$. Then we may
proceed to define ${\cal J}_1 = (0,\tau^{-1} R')$, ${\cal J}_p$ as the
preimage of ${\cal J}_1$ by $H$ inside $J_p$, and for $1<q<p$, the
interval ${\cal J}_q$ is the
preimage of ${\cal J}_p$ by $H^{p-q}$ inside $J_q$. Since we decreased
the intervals, ${\cal J}_q$ are pairwise disjoint and contained in
$(0,R')$.

Let us now define the ``presentation function'' $\Pi$, initially only on the
union of intervals ${\cal J}_q$. We put $\Pi(x) = \tau x$ for $x\in
{\cal J}_1$, $\Pi(x) = H(x)$ if $x\in {\cal J}_q$, $1<q<p$, and
$\Pi(x) = \tau H(x)$ if $x\in {\cal J}_p$. We use notation: $\Pi_q$ is
the restriction of $\Pi$ on ${\cal J}_q$, $1\le q\le p$.

\paragraph{The analytic continuations of $\Pi$, $\hat\Pi$.}
To define the analytic continuation of $\Pi$ more precisely, consider
the following geometrical disks: $D_1 = D(0, R')$, $D_2=\tau^{-1}D_1$.
As in the preceding paragraph, $R'$ is less than $b_0$ but large
enough so that $J_i \subset [0, R']$ for $i=1,\cdots,p$.  
Then $\hat{D}_i$ are
analogously defined disks in the phase space of $\hat{H}$. Now we
consider the analytic continuation of $H$ to the following sets. We
extend the linear branch $\Pi_1$ to $U_1 := D_2$. Then $\Pi_p$ is
extended to the ``figure eight'' set $U_p$ chosen so that 
$H$ restricted to each connected component of $U_p$ is a covering of
the punctured disk centered at $0$ with radius $\tau^{-1} R'$.  
From the limit formula for $H$ in
Definition~\ref{defi:10fa,1}, $U_p$ is contained in the geometric disk
with diameter ${\cal J}_p$. Then for $1<q<p$ we set $U_q =
H^{q-p}(U_p)$ choosing the appropriate connected component of the
preimage, which contains the interval ${\cal J}_q$. 

Observe that domains $U_i$ do not intersect.
By hypotheses of the EW-class, see Definition~\ref{defi:10fa,1}, 
sets $U_2, \cdots, U_p$
are contained in geometric disks based on the corresponding ${\cal
J}_q$, and so are pairwise disjoint and also disjoint with $U_1 =
D_2$.

\paragraph{Plan of the proof.}
$\Pi$ is defined by analytic continuation to $\bigcup_{q=1}^{p} U_q$ and     
the same construction can be carried out for $\hat{H}$, yielding
a box mapping $\hat{\Pi}$. 

See Figure~\ref{fig:16fa,1} for an illustration in the case of $p=3$. 

\begin{figure}
\epsfig{figure=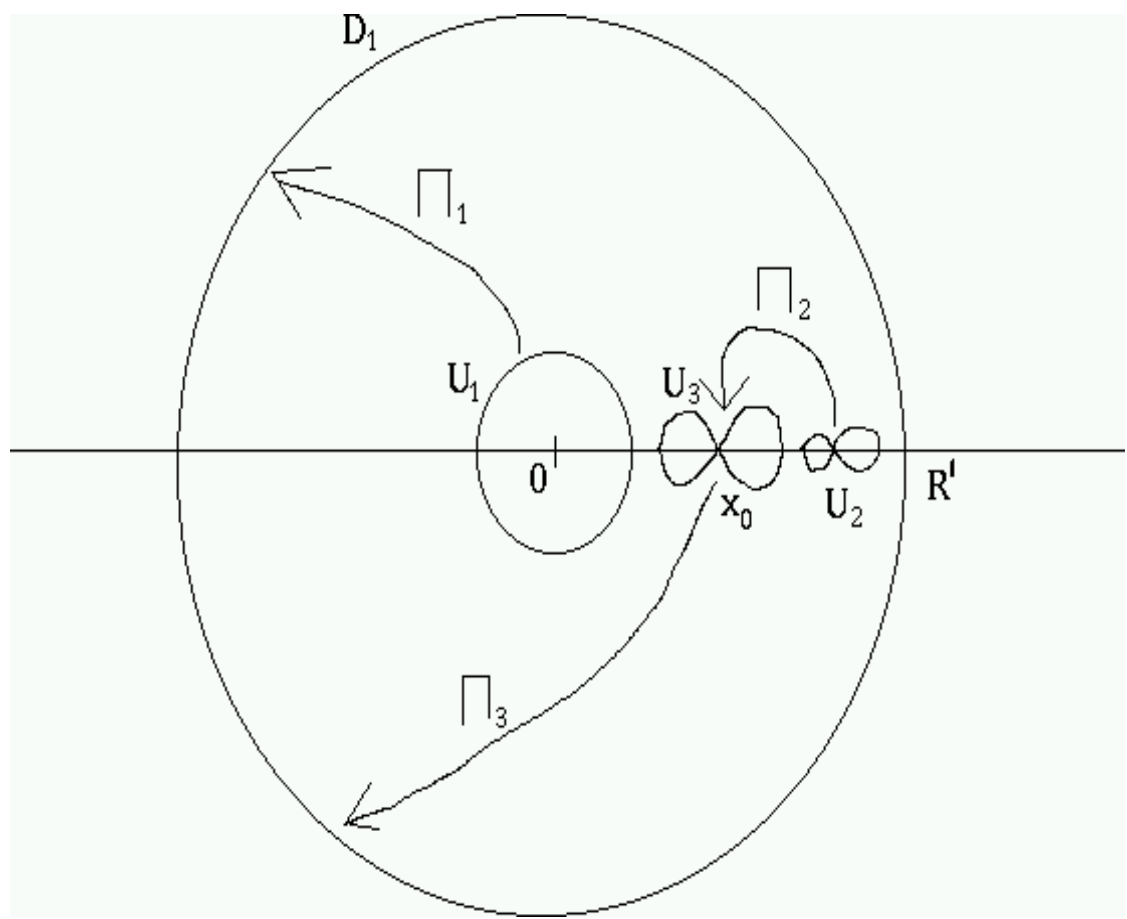}
\caption{The box mapping $\Pi$.\protect\label{fig:16fa,1}}
\end{figure}

The partial conjugacy referred to by
Proposition~\ref{prop:11fp,2} is then obtained as the conjugacy
between the box mappings $\Pi$ and $\hat{\Pi}$. One might wonder how
that is possible, since analytically $\Pi$ is no simpler than $H$,
having exactly the same type of singularity at $x_0$. The answer is
that the dynamics of $\Pi$ is completely different from $H$. In
particular, $0$ has become a repelling fixed point, and so $\Pi$ is a
post-critically finite map, making the task of constructing the
conjugacy much easier, again using the pull-back method.

\paragraph{Preparatory estimates.}
\begin{lem}\label{lem:16fa,1}
Suppose that $g$ is real-analytic at $0$ with the following
power-series expansion: 
\[ g(x) = x - \varepsilon x^3 + O(|x|^4)\; ,\]
with $\varepsilon>0$. 

If $a_2 < 0 < a_1$ are in the basin of attraction of $0$, then there
exists $K>0$ such that for every $n\geq 0$ 

\[ K \leq \frac{|g^n(a_2)|}{|g^n(a_1)|} \leq K^{-1} \; .\]
\end{lem}
\begin{proof}
The Fatou coordinates $h_{\pm}$ on the right and left attracting petal
of $g$, respectively, are $\frac{1}{2\varepsilon z^2} + O(|\log
z|)$ with the leading term the same on either side of $0$,
see~\cite{carleson}. It follows that for $n$ sufficiently large 
\[ \frac{1}{2} \sqrt{\frac{n}{2\varepsilon}} < |g^n(a_i)| <
2\sqrt{\frac{n}{2\varepsilon}} \]
and the lemma follows.
\end{proof}

\begin{lem}\label{lem:16fa,2}
Suppose now that two mappings $g, \hat{g}$ are given, both analytic in
a neighborhood of $0$, $g$ in the same form as in
Lemma~\ref{lem:16fa,1} and $\hat{g}$ in the analogous form:
\[ \hat{g}(x) = x - \hat{\varepsilon} x^3 + O(|x|^4) \; ,\]
$\hat{\varepsilon} > 0$. Let $h_{\pm}, \hat{h}_{\pm}$, respectively,
denote ``right'' and ``left'' Fatou coordinates symmetric about the
real axis. Suppose that $\Upsilon$
maps $0$ to itself, is $\hat{h_+}^{-1} \circ h_+$ to the right of
$0$ and $\hat{h_-}^{-1}\circ h_-$ to the left. Then $\Upsilon$ is
quasi-symmetric in a neighborhood of $0$. 
\end{lem}
\begin{proof}
We observe first that $\hat{h_+}^{-1}\circ h_+$ and $\hat{h_-}^{-1} \circ h_-$
are quasisymmetric in the respective one-sided neighborhoods of $0$. 
This follows from the fact proved in~\cite{carleson} that each of the Fatou
coordinates has the form $\Gamma(z^2)$ with $\Gamma$ quasi-conformal
on the plane, which can be normalized to a map from the left real semi-line
into itself. It remains to show, see~\cite{jaswi} Lemma 3.14,  
that for all $0<\alpha<\alpha_0$, with
$\alpha_0$ chosen conveniently small, and fixed $K>0$
\begin{equation}\label{equ:16fa,1} 
K^{-1} < \frac{|\Upsilon(\alpha)|}{|\Upsilon(-\alpha)|} < K \; .
\end{equation}

Because of the symmetry between the right and left side, we will only
show the lower estimate. 
To this end, fix some $a_1 > 0$ and $a_2 < 0$ and set $\hat{a}_i =
\Upsilon(a_i)$, $i=1,2$. Without loss of generality $|a_1|, |a_2| >
\alpha_0$. Find the smallest $n$ such that $g^n(a_1) \leq
\alpha$. Then $g^n(a_1)/\alpha > 1/2$ if $\alpha_0$ was small enough
so $|g^n(a_2)| > K_1 \alpha$ for some fixed $K_2 > 0$ based on
Lemma~\ref{lem:16fa,1}. Since $\Upsilon(g^n(a_2)) =
\hat{g}^n(\hat{a}_2)$ and the left branch of $\Upsilon$ is
quasisymmetric, we get 
\[ |\hat{g}^n(\hat{a}_2)| > K_2 |\Upsilon(-\alpha)| \]
with fixed $K_2 > 0$. But finally 
\[ \Upsilon(\alpha) \geq \hat{g}^n(\hat{a}_1) \geq K_3
|\hat{g}^n(\hat{a}_2)| \]
with $K_3 > 0$ depending only on the choice of $a_1, a_2$ by
Lemma~\ref{lem:16fa,1}. The lower estimate of
inequality~(\ref{equ:16fa,1}) follows. 
\end{proof}

\paragraph{Construction of the partial conjugacy.}
We will now resume work on proving Proposition~\ref{prop:11fp,2} first
by building a partial conjugacy between $\Pi$ and $\hat\Pi$. 
We start by considering the affine map $\varphi_0(z) =
\frac{\hat{R'}}{R'}z$.  

Next, we will construct a quasiconformal map $\varphi_1$. It is not
going to be defined on the entire plane. Outside of $D_1$, 
we set $\varphi_1 = \varphi_0$. On $U_1$, $\varphi_1(z) =
\frac{\hat{R'}}{R'} \frac{\tau}{\hat{\tau}}z$. This will ensure 
$\varphi_1 \circ \Pi_1 = \hat{\Pi}_1 \circ \varphi_1$ on $U_1$.  
Then on $U_p$ we make
$\varphi_1$ equal to the lift of the affine $\varphi_{1|_{U_1}}$ 
by $H$, $\hat{H}$,
set up so that the lifted mapping sends $U_p \cap \RR$ into $\RR$ 
preserving the orientation. Equivalently, this is the lifting of
$\varphi_{0|D_1}$ by $\Pi_p, \hat{\Pi}_p$. 
Finally, on each $U_q$, $1<q<p$ we set 
$\varphi_1 = \hat{H}^{q-p} \varphi_1 H^{p-q}$ applying the appropriate
inverse branch. Summarizing, $\varphi_1$ is symmetric about the real
line, fixes $0$, inside $D_2$ is defined on
the union of sets $U_1, \cdots, U_p$ and satisfies $\varphi_1 \Pi =
\hat{\Pi} \varphi_1$ on the union of their boundaries. What we still
need is extend the domain of definition of $\varphi_1$ to the entire
plane. 

Before we do, observe that $\varphi_1$ restricted to the real line is
quasi-symmetric provided that we interpolate on the intervals where it
has not been defined, for example, by affine maps. This is clear, since 
on $D_1 \cap \RR$ the map $\varphi_1$ is piecewise analytic and at the
point of contact of two pieces usually it can be continued from either
of them to a neighborhood of its closure. 
An exception occurs if the common endpoint is 
$x_0$ or one of its preimages. However, in the
neighborhood of $x_0$ we can invoke Lemma~\ref{lem:16fa,2}, and the
map has been propagated to the preimages of $x_0$ by diffeomorphic
branches of $H, \hat{H}$. This allows us to construct a quasiconformal
homeomorphism $\varphi_2$, of the lower half-plane onto itself, whose
continuous extension matches $\varphi_1$ on the real line. 

\begin{figure}
\epsfig{figure=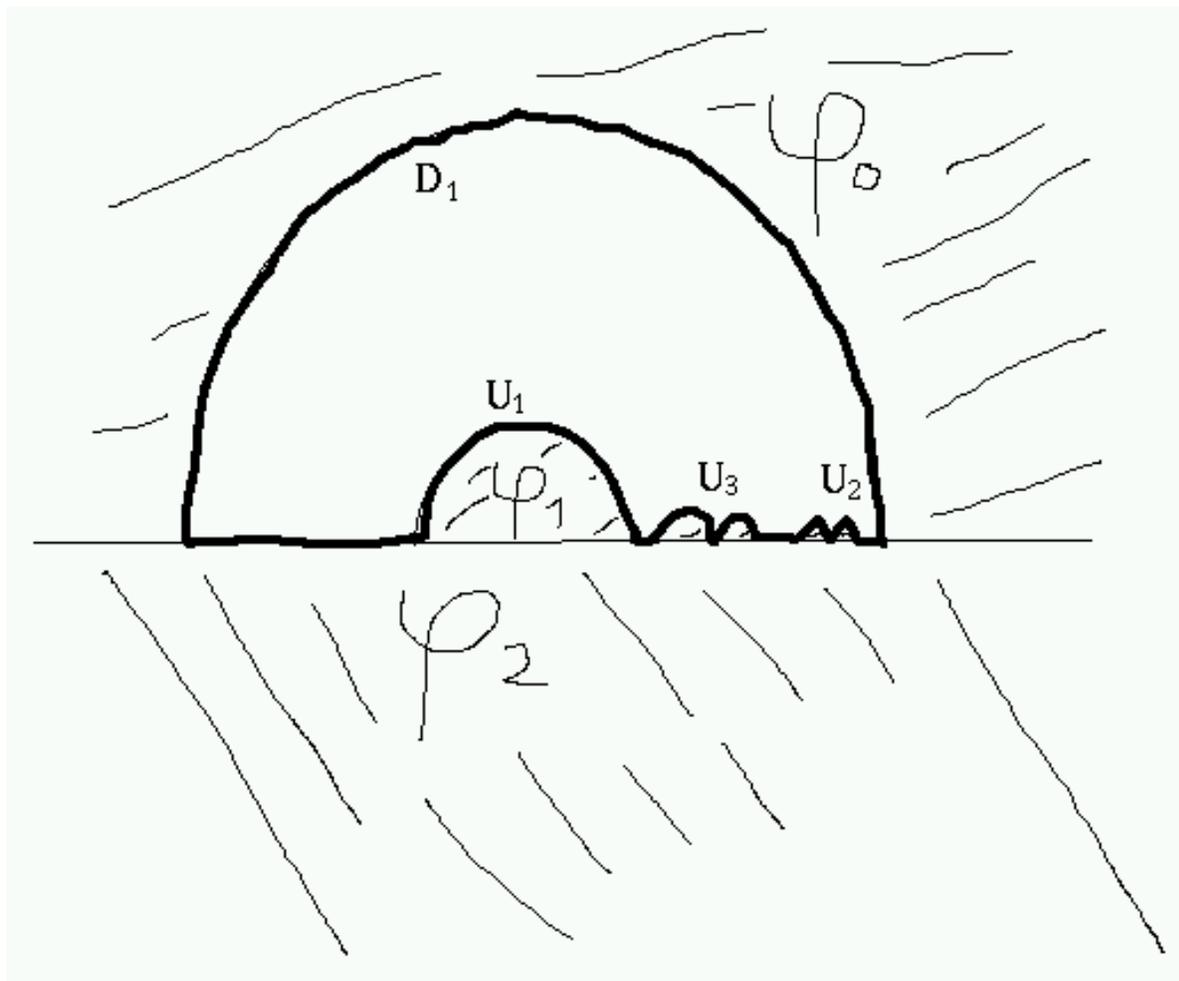}
\caption{Construction of the partial conjugacy.\protect\label{fig:23fa,1}}
\end{figure}

Now the reader is invited to consult Figure~\ref{fig:23fa,1} and pay
attention to the curve $w$ marked by a thick line. This line 
consists of the boundary curves of domains $U_1, \cdots, U_p$
intersected with the upper half-plane, pieces of the real line between
them and the boundary arcs of $D_1$. The key
fact about $w$ is that it is a quasi-circle. Indeed, it consists of
finitely many quasi-conformally embedded arcs intersecting
always with a certain angle fitting between them. In particular, at
$x_0$ the curves are still known to posses tangent lines making angles
$\pi/4$ with the real line, see~\cite{carleson}. Similarly, the curve 
$\hat{w}$ built of the analogous arcs in the phase space of $\hat{H}$,
with the short-cut which is the image of the corresponding part of $w$
by $\varphi_0$, is also a quasi-circle. 

Next, we define a quasi-conformal map $\varphi_3$ of the unbounded component of
the complement of $w$ onto the unbounded component of the complement
of $\hat{w}$. In the lower half-plane, we set $\varphi_3 =
\varphi_2$. On $\HH^+ \setminus D_1$ we set $\varphi_3 =
\varphi_0$. 
On $U_i \cap \HH^+$, $i=1,\cdots,p$, 
we make $\varphi_3 =
\varphi_1$. Note that we get that 
$\varphi_3$ extends $\varphi_{1|\HH^+}$ . But now
$\varphi_3$ can be extended to the entire plane by reflecting about the
quasi-circles $w, \hat{w}$. Finally, we take $\varphi_3$ from $\HH^+$
and reflect it about the real line to $\HH^-$, thus obtaining the
desired quasi-conformal extension of the map $\varphi_1$ to the entire
plane. We will still use the notation $\varphi_1$ for this extension. 

\paragraph{Pull-back.}
Now we use the standard pull-back construction. That is, we construct
a sequence of quasi-conformal homeomorphisms of the plane
$\Upsilon_n$, $n=0,1,\cdots$, with $\Upsilon_0 = \varphi_1$ and
$\Upsilon_n$ re-defined on each $U_i$, $i=1,\cdots, p$ according to
the formula $\hat{\Pi}^{-1} \circ \Upsilon_{n-1} \circ \Pi$. This is
rigorous except for $i=p$ where $H^{-1}, \hat{H}^{-1}$ cannot be
defined and one should talk instead of the lifting of $\Upsilon_{n-1}$
to universal covers. Since $\varphi_1(0) = 0$ and that condition is
preserved by the pull-back, the lifting is well defined.    

$(\Upsilon_n)$ form a compact family of homeomorphisms, so we can find
$\Upsilon_{\infty}$ which is the limit of some subsequence of them. Then
$\Upsilon_{\infty}$ is quasi-conformal and still coincides with
$\varphi_1$ on $\partial U_i$. 

\paragraph{Dynamics of $\Pi$.}
We want to show that $\Upsilon_{\infty}(H^n(1)) = \hat{H}^n(1)$ for
all non-negative $n$. Since $\Upsilon_{\infty}(0)=0$ as well, this
will mean that $\Upsilon_{\infty}$ conjugates the forward critical
orbits. The dynamics on the critical orbit under $\Pi$ is simple to
understand using the functional equation: $1$ is periodic with period
$p$ and every image of $1$ is eventually mapped to $1$. Let us
consider the filled-in Julia $K_{\Pi}$ defined as the set of all points
which can be forever iterated by $\Pi$. Every point $x\in K_{\Pi}$ has
an itinerary consisting of symbols $1,\cdots,p$, where the $k$-th
symbol being $i$ means that $\Pi^k(x) \in U_i$. The key observation is
that no two points can have the same itinerary. This follows because
$\Pi$ expands, though not uniformly, the hyperbolic metric of the
punctured disk $V := D(0,R') \setminus \{0\}$. Indeed, the map $H^{p-q+1}$
for $q>1$ is a covering of $V$ by $U_q$. Only on $U_1$ is $\Pi$ an
isometry. But every point is $K_{\Pi}$ can only be iterated by $\Pi_1$
finitely many times, and it follows that the expansion ratio with
respect to the hyperbolic metric along the orbit of any $x\in K_{\Pi}$
goes to $\infty$. Since all $U_q$, $q>1$, have finite diameters with
respect to this metric, it follows that the distance between two
points with the same itinerary must be $0$.  

The partial conjugacy $\Upsilon_n$ preserves the first $n$ symbols of
any itinerary. So $\Upsilon_{\infty}$ maps $K_{\Pi}$ into
$K_{\hat\Pi}$ preserving the itineraries. But $H^n(1)$ and
$\hat{H}^n(1)$ have the same itineraries, so it follows that
$\Upsilon_{\infty}(H^n(1)) = \hat{H}^n(1)$.

\paragraph{Use of the functional equation to finish the proof.}
Mapping $\Upsilon_\infty$ satisfies the dynamical condition required on
$\phi_1$ in the statement of Proposition~\ref{prop:11fp,2}, but has no
reason to obey the requirement imposed on the boundary of $U_- \cup
U_+$. To correct this, first restrict $\Upsilon_{\infty}$ to the set
$U_p$. Such a restriction still satisfies $\Upsilon_{\infty}(H^p(z)) = 
\hat{H}^p(\Upsilon_{\infty}(z))$ for every $z$ in the forward orbit
of $x_0$ by $H^p$. Additionally, by our construction, it also satisfies
$\frac{\hat{R'}\tau}{\hat{\tau}R'} H(z) = \Upsilon_{\infty}(\hat{H}(z))$ on the
boundary of $U_p$, which gets mapped on the geometric circle $C(0,
\tau^{-1} R')$ by each branch of $H$. On the annulus $\{ z :\: \tau^{-1}R' \leq |z| \leq \tau^{-1}R\}$
we can define a quasi-conformal map $\upsilon$ which is linear with
slope $\frac{\hat{R'}\tau}{\hat{\tau}R'}$ on $C(0, \tau^{-1}R')$ and linear with slope
$\frac{\hat{R}\tau}{\hat{\tau}R}$ on $C(0,\tau^{-1}R)$. 
Taking the appropriate lift
$\hat{H}^{-1} \circ \upsilon \circ H$, we can modify
$\Upsilon_{\infty}$ to a new map $\phi'_1$ which is defined on
$U' := H^{-1}(D(0,\tau^{-1}R))$, is the same as $\Upsilon_{\infty}$, in
particular conjugating forward critical orbits of $x_0$ by $H^p,
\hat{H}^p$, on $U_p$ and satisfies $\frac{\hat{R}\tau}{\hat{\tau} R}
H(z) = \hat{H}(\phi'_1(z))$ on the boundary of $U'$. 

Note that $\phi'_1$ restricted to $U'_1\cap \RR$ is quasi-symmetric. Indeed,
$\phi'_1$ restricted to a smaller interval $U_p \cap \RR$ was just a
restriction of a quasi-conformal homeomorphism of the plane. We then
extended it to a larger interval $U'$ and the new mapping
remains quasi-symmetric since it extends quasi-conformally to a
neighborhood of each of the endpoints of $U_p$.

Finally, $\phi_1$ as postulated by Proposition~\ref{prop:11fp,2} is
given by the formula 
\[ \phi_1 = \hat{G}^{-1} \circ \phi'_1 \circ G \; .\]
Immediately, we see that $\phi_1$ is quasi-symmetric when restricted
to $U \cap \RR$ since it is just
the pulled-back of a quasi-symmetric mapping from
$U'$ by analytic maps $G, \hat{G}$. 

We check the conditions starting from the functional equation
$\tau^{-1} H = H\circ G$ satisfied on $U := U_-\cup U_+$. First, 
$G^{-1}(U') = (H\circ G)^{-1}(D(0, \tau^{-1}R)) = H^{-1}(D(0,R)) =
U$ so the domain of $\phi'$ is $\overline{U}$ and, by an analogous
argument, its range is $\overline{\hat{U}}$. For $z\in \partial U$, 
\[ \hat{H}(\phi_1(z)) = \hat{H} \circ \hat{G}^{-1}(\phi'_1 (G(z))) =
\hat{\tau} \hat{H} \circ \phi'_1(G(z)) = \hat{\tau} \frac{{\hat
R}\tau}{\hat{\tau} R} H(G(z)) = \frac{\hat{R}}{R} H(z) \]
as needed. 

To verify the conjugacy on the forward critical orbit, we use the
identity $H^p \circ G = G \circ H$ valid at least on $[0, b_0]$, see
the beginning of the proof of Proposition~\ref{prop:10fp,1}. Thus, 
\[ \phi_1 H^n(x_0) = \hat{G}^{-1} \circ \phi'_1 ( G(H^n(x_0)) =
\hat{G}^{-1} \circ \phi'_1(H^{pn}(G(x_0))) = \]
\[ =\hat{G}^{-1}\circ
\hat{H}^{pn}(\phi'_1(G(x_0))) = \hat{H}^n(\hat{G}^{-1}(\phi'_1(G(x_0)))) = 
\hat{H}^n \phi_1(x_0) \; .\]

This concludes the proof of Proposition~\ref{prop:11fp,2}.

\paragraph{Extension of $\phi_1$.}
The derivation of Proposition~\ref{prop:11fp,1}
from Proposition~\ref{prop:11fp,2} is another standard
application of the pull-back method. First, however, we have to extend
$\phi_1$ obtained from Proposition~\ref{prop:11fp,2} to the complex
plane in such a way as to make a conjugacy on the boundary of $U_-
\cup U_+$. In view of the claim of Proposition~\ref{prop:11fp,2}, we
simply need to extend $\phi_1$ to the whole plane in such way that it
becomes linear with slope $\frac{\hat{R}}{R}$ outside of $D(0,R)$, so
that the main difficulty is interpolating on $D(0,R) \setminus U$. 

First, we perform this interpolation on the real line, constructing a
quasi-symmetric homeomorphism $\varphi_1$ which coincides with
$\phi_1$ on $U\cap\RR$ and is linear with slope $\frac{\hat{R}}{R}$
outside $(-R, R)$. Next, we extend $\varphi_1$ quasi-conformally to
the lower half-plane, getting a
picture similar to one shown of Figure~\ref{fig:23fa,1}. By now, we
have a quasi-conformal map defined on the complement of $\HH^{+}
\setminus U$. But the boundary of $\HH^+ \setminus U$ is a
quasi-circle, for the same reasons as the curve $w$ in the proof 
of Proposition~\ref{prop:11fp,2}. So we can extend this to a homeomorphism 
of the plane by quasi-conformal reflection. Finally, we make the
mapping symmetric about the real axis by reflecting from the upper
half-plane into the lower. This gives the extension of $\phi_1$ with
the desired properties: it is a quasiconformal homeomorphism of the
plane and the conjugacy condition $\phi_1 (H(z)) = \hat{H}(\phi_1(z))$
now holds on the boundary of $U$ as well as on the forward orbit of
$x_0$.

\paragraph{Proof of Proposition~\ref{prop:11fp,1}.}
Thus, we construct a sequence of
quasi-conformal homeomorphisms $\phi^n$ of the plane, by setting
$\phi^0 = \phi_1$ and defining $\phi^n$ for $n>0$ as $\phi^{n-1}$
outside of $U_+\cup U_-$ and to be the lifting of $\phi^{n-1}$ to the
universal covers $H_{|U_+}, \hat{H}_{|\hat{U}_+}$ and $H_{|U_-},
\hat{H}_{|\hat{U}_-}$. Both the lifting are uniquely defined by the
requirement that $\phi^n$ should fix the real line with its
orientation. 

The sequence $\phi^n(z)$ actually stabilizes for every $z\notin
K_H$. So $\phi^n$ converge on the complement of $K_H$ and by taking a
subsequence can be made to converge globally to some map
$\phi^{\infty}$. Outside of $K_H$, $\phi^{\infty}$ satisfies the
functional equation $\phi^{\infty} H = \hat{H} \phi^{\infty}$ and then
it also satisfies it on $K_H$ by continuity, in the light of
Theorem~\ref{pre}. So we can set $\phi_0 := \phi^{\infty}$ and this
concludes the proof of Proposition~\ref{prop:11fp,1}.

\section{Rigidity}
In this section we will prove Theorem~\ref{theo:10fa,2} by
constructing towers based on two EW-maps and showing that they must be
the same. 

\subsection{Towers and their dynamics} 
Let $H$ belong to the EW-class with some combinatorial type $\aleph$.

\begin{defi}\label{defi:11fa,1}
Define, for $n=0,1,2,...,$
$H_n(z)=\tau^n H(z/\tau^n)$.
Then $\tau^n K_H$ is the Julia set of the map
$H_n:U_n\to V_n$, where $U_n=\tau^n(U^+\cup U^-), V_n=\tau^n V$. 
Note that, for any $n>m$,
$H_m=H_n^{|\aleph|^{n-m}}$.

The collection of maps $H_n:U_n\to V_n$, $n=0,1,...$
forms the {\it tower} of $H$.
\end{defi}

It is important to realize that $H_{n+1}^{|\aleph|} = H_n$ for all $n=0,
1,\cdots$. Each $H_n$ has its filled-in Julia set $K_{H_n}$, see
Definition~\ref{defi:29fa,1}. 
It follows straight from the definition of $H_n$, that
$K_{H_n} = \tau^n K_{H}$. 
Another property which follows from
the definition is that the sequence $K_{H_n}$ is increasing with $n$. In
line with the general strategy of working with towers, we will need this: 

\begin{prop}\label{prop:11fa,3}
In the tower of every EW-map $H$, the Julia set
\[ \bigcup_{n=1}^{\infty} \tau^n K_{H} \]
is dense in $\CC$. 
\end{prop}

\paragraph{Dynamics in towers.}
{\em Tower dynamics} is understood as the set of all possible
compositions of mappings $H_i$ from the tower. So, if we say that $z$
is mapped to $z'$ by the tower dynamics, it means that a composition
exists which sends $z$ to $z'$.  
The key statement about the dynamics in towers generalizes
Lemma~\ref{lem:28fp,2} and uses the same notation.

Introduce the following sets.
Let $\omega_n$ be the omega-limit set of $0$ under the action of $H_n$.
In particular, $\omega_0=\omega$. Each $\omega_n$ is a closed set.
Introduce $\omega_\infty=\cup_{n\ge 0}\omega_n$.
It is also a closed subset of the plane.
Furthermore, $\omega_\infty\cap V_n
=\omega_\infty\cap U_n=\omega_n$.

\begin{prop}\label{prop:11fa,1}
For every  $z\in \CC$ which is never mapped to $\RR$ by the tower
dynamics, there exist
sequences $z_n \in \CC$ and $m_n \in \NN\cup \{0\}$, $n=0,1,\cdots$,
such that $z_0 = z$, $z_n$ is an image of
$H_{m_{n-1}}(z_{n-1})$ by the tower dynamics, for every $n>0$, and at least
one of the following statements is true: 
\begin{itemize}
\item
there exists $\eta>0$ such that 
$\dist (z_n, \omega_{\infty}) > \eta \tau^{m_n}$ for every
$n>0$, with $\dist$ meaning the Euclidean distance, or
\item
for every $n>0$ 
\[ \tau^{-m_n} z_n \in (U_-\cup U_+) \setminus
(U_{+,c} \cup U_{-,c})\; .\]
\end{itemize}
\end{prop}

To prove that one of the alternative statements must hold, notice
first that without loss of generality $z\notin K_{h}$ for any
$h$. Otherwise, the alternative will follow by applying
Lemma~\ref{lem:28fp,2} inductively to the dynamics $H_h$. 

So, assuming that $z_{n-1}$ has been constructed we map it by the
dynamics of $H_{m_{n-1}}$ until the first moment $q$ when $w :=
H^q_{m_{n-1}}(z)$ is no longer in the domain of $H_{m_{n-1}}$. 
The only point where the set $D(0,R) \setminus (U_-\cup U_+)$ touches
$\omega_{\infty}$ is $x_0$. So, if $w \notin \tau^{m_{n-1}} D(x_0,
\varepsilon)$ for some $\varepsilon>0$, 
then $w \in \tau^{m_{n-1}+m_0} (U_- \cup U_+)$ and
$\dist(w,\omega_{\infty}) > \tau^{m_{n-1}} \eta$ with $m_0$ and
$\eta>0$ which depend only on $\varepsilon$. In that case we set $z_n
:= w$ and $m_n = m_{n-1} + m_0$.

Otherwise, we continue iterating $W := \tau^{-m_{n-1}}w$ by $G$. 
The connection with the
tower dynamics  relies on the following
simple observation:

\begin{fact}\label{fa:5gp,1}
For any $q, Q$, the composition $G^q(\tau^{-Q}z)$ can be represented
as $\tau^{-s} \chi(z)$ where $\chi$ belongs to tower dynamics.
\end{fact}
\begin{proof} 
If $G^{q-1}(\tau^{-Q} z) = \tau^{-s'} \chi'(z)$, then
\[ G^q(\tau^{-Q}z) = G (G^{q-1}(\tau^{-Q}z)) = H^{p-1}(\tau^{-s'-1}
\chi'(z)) = \tau^{-s'-1} H_{s'+1}^{p-1}(\chi'(z)) \; .\]
\end{proof}

We will continue iteration by $G$ until the first moment
$q'$ when $W' := G^{q'}(W)$ is either outside of
$D(x_0, \varepsilon)$, or the distance from $\arg (W' - x_0)$ 
to $0$ or $\pi$ on the circle is less than
$\pi/5$.    

By specifying $\varepsilon$ to be sufficiently small, we can achieve
the following for every $u\in D(x_0, \varepsilon)$, $u\neq x_0$:
\begin{itemize}
\item
$|\arg (G(u) - x_0) - \arg(u-x_0)| < \pi/10$,
\item
if the distance from $\arg(u-x_0)$ to $0$ and $\pi$ on the circle is
less than $\pi/5$, then $u\in U_- \cup U_+$  
\item
if $u \in U_{-,c} \cup U_{+,c}$, then the distance from $\arg(u-x_0)$ to
$0$ or $\pi$ is less than $\pi/10$, 
\item $G(u)$ is in $\tau U_-$. 
\end{itemize}

The first possibility is that $|W'-x_0| \geq \varepsilon$. By the
properties postulated here, the distance from $\arg(W' - x_0)$ to $0$
and $\pi$ on the circle is at least $\pi/10$ and $W' \in \tau U_-$. By
Fact~\ref{fa:5gp,1}, for some $s$ we get $z_n := \tau^s W' = \chi'(w)$      
for some tower iterate $\chi'$. Then $z_n \in \tau^{s+1} U_-$ and 
$\dist(z_n, \omega_{\infty}) \geq
\tau^{s}\varepsilon\sin\frac{\pi}{10}$. We set $m_n = s + 1$. 

Finally, it may be that $|W'-x_0| < \varepsilon$. Then the distance
from $\arg(W' - x_0)$ to $\{0,\pi\}$ on the circle is between $\pi/10$
and $\pi/5$. By the choice of $\varepsilon$, $W' \in (U_-\cup U_+)
\setminus (U_{+,c} \cup U_{-,c})$. Again, we set $z_n = \tau^s W'$
where $s$ comes from Fact~\ref{fa:5gp,1}. Setting $m_n = s$, we get
\begin{equation}\label{equ:5gp,1}
\tau^{-m_n} z_n \in (U_-\cup U_+) \setminus
(U_{+,c} \cup U_{-,c}) \; .
\end{equation}

This inductive construction yields a sequence of points $z_n$ and
integers $m_n$ such that for each of them either $\dist(\tau^{-m_n}
z_n, x_0) > \eta$ with $\eta$ independent of $n$, as happens in the
first two cases we considered, or $z_n$ satisfies
condition~(\ref{equ:5gp,1}). 
Since one of these subsequences is infinite,
Proposition~\ref{prop:11fa,1} follows.

\subsection{Expansion of the hyperbolic metric.}
\paragraph{Hyperbolic metric.}

Recall that $\omega_n$ is the omega-limit set of $0$ under the action of $H_n$,
$\omega=\omega_0$, and $\omega_\infty=\cup_{n\ge 0}\omega_n$.

Let $\rho_\infty$
be the hyperbolic metric of $\CC\setminus \omega_\infty$.
Note that $\rho_\infty$ is invariant under the rescaling
$z\mapsto \tau z$. 

The following lemma is stated in terms of $H$, but clearly it applies
to any $H_k$ as well, because the only difference is the conjugation
by a power of $\tau$, which is the isometry of the hyperbolic metrics
involved. 

\begin{lem}\label{lem:5ga,1} 
Suppose that $H$ is an EW-map with combinatorial type $\aleph$. 
For any $z\in (U_-\cup U_+) \setminus H^{-1}(\omega_{\infty})$, we get 
that the hyperbolic metric expansion ratio 
\[ DH_{\rho_{\infty}}(z) \geq (\iota'(z))^{-1} \]
where $\iota$ is the inclusion map from $\CC \setminus
H^{-1}(\omega_{\infty})$ into $\CC \setminus \omega_{\infty}$ and the
prime denotes its contraction ratio with respect to the corresponding
hyperbolic metrics. 
\end{lem}
\begin{proof}
We can represent $H'(z) = DH_{\rho_{\infty}}(z)  \iota'(z)$
where $H'(z)$ represents the expansion ratio of $H$ acting from the
hyperbolic metric of $\CC \setminus H^{-1}(\omega_{\infty})$ into the
hyperbolic metric of $\CC \setminus \omega_{\infty}$. Writing
$p:=|\aleph|$, we get for any $k\geq
0$ that $H = H_k^{p^k}$. Observe that $H_k^{p^k}$ 
is a holomorphic covering of 
$X_k = D(0, \tau^k R) \setminus \omega_{\infty}$ by 
$\tau^k (U_-\cup U_+)\setminus
H_k^{-p^k}(\omega_{\infty})$.

Hence, it is a local isometry with respect
to the corresponding hyperbolic metrics. So, it is non-contracting
when the hyperbolic metric of 
$\tau^k (U_-\cup U_+) \setminus H_k^{-p^k}
(\omega_{\infty})$ is replaced with the hyperbolic metric of a larger
set 
$Y_k = \tau^k (U_-\cup U_+)\setminus H^{-1}(\omega_{\infty}))$.

As $k$
tends to $\infty$, the hyperbolic metrics of $X_k$ tend to
$d\rho_{\infty}$ while the hyperbolic metrics of $Y_k$ tend to the
hyperbolic metric of $\CC \setminus H^{-1}(\omega_{\infty})$ uniformly
on compact sets. It follows that $H'(z) \geq 1$ as needed. 
\end{proof}

\paragraph{Uniform expansion.} 
Now take any point $z\in \CC$ which is never mapped to $\RR$ by the
tower dynamics. Proposition~\ref{prop:11fa,1} then delivers a sequence
$z_n$. Let $\chi_n$ be the corresponding tower iterate which maps $z$
to $z_n$.

\begin{lem}\label{lem:5gp,1}
For every $D$ there exists $\lambda >
1$, such that for every $n$ and every $w$ in the ball centered at $z_n$
with radius $D$ with respect to $\rho_{\infty}$,
$D_{\rho_{\infty}}H_{m_n}(w) > \lambda$, provided that $w$ is in the
domain of $H_{m_n}$. 
\end{lem}
\begin{proof}
By Lemma~\ref{lem:5ga,1} and
Fact~\ref{fa:20gp,2}, $DH_{\rho}(w) > \lambda > 1$ where $\lambda$
depends only on the distance in $\rho_{\infty}$ from $w$ to
$H^{-1}(\omega_\infty)$. By rescaling, the same is true for all $H_k$. But if
either case of the alternative statement holds, points $z_n$ are
all in a uniformly bounded $\rho_{\infty}$-distance from the
corresponding set $H_{m_n}^{-1}(\omega_{\infty})$. The same will be true 
for $w$ by the triangle inequality. 
\end{proof}

\begin{lem}\label{lem:5gp,2}
For every $n$, let $\zeta_n$ denote the inverse branch of $\chi_n$ which
maps $z_n$ to $z$ defined on some simply-connected set $U_n \ni z_n$. 
 Then for every $D$ and $\varepsilon$ there exists
$n_0$ such that for every $n\geq n_0$ if the diameter of $U_n$ with
respect to $\rho_{\infty}$ does not exceed $D$, then $\zeta_n(U_n)$ is
inside the hyperbolic ball of radius $\varepsilon$ centered at $z$. 
\end{lem}
\begin{proof}
Pulling back a $U_n$ will not increase its
diameter, so each time we pass $z_m$ its radius will be shrunk by a
definite factor.
\end{proof}

\paragraph{Density of the Julia sets.}
We can now prove Proposition~\ref{prop:11fa,3}. For some fixed $D$ and
every $n$, we can find 
an element of $H^{-1}_{m_n}(\omega_{\infty})$, moreover, a preimage of
$0$ by $H_{m_n}$,  which can be joined 
to $z_n$ by a simple arc $\gamma_n$ of hyperbolic length which does not exceed 
some fixed $D$ and which is completely contained in $\tau^{m_n}
(U_-\cup U_+)$. This follows from simple geometric considerations
similar to those used in the proof of Theorem~\ref{pre}. 
We can then find $k$ which is at least equal to $m_n$
and large enough so that the tower iterate $\chi_n$ can be represented
as an iterate of $H_k$. 

Then the inverse branch $\zeta_n$ is defined on a neighborhood of
$\gamma_n$. We can apply Lemma~\ref{lem:5gp,2}  to get that $\zeta_n$ maps
$\gamma_n$ into a neighborhood of $z$ whose diameter shrinks to $0$ as
$n$ grows. Letting $n$ go to
$\infty$, we get that every ball
centered $z$ contains a preimage of $0$ by some iterate of the tower dynamics.
But every preimage of $0$ in the tower belongs to some
$K_{H_k}$ and so Proposition~\ref{prop:11fa,3} follows.

\subsection{Conjugacy between towers}
Given towers built for two EW-maps $H$ and $\hat{H}$, we construct a
quasiconformal conjugacy {\em between the towers}  by rescaling  the
conjugacy between $H$ and $\hat{H}$  to 
conjugacies $\tau^n\circ \phi_0\circ \tau^{-n}$
of $H_n, \hat H_n$, pass to a limit, and 
get a conjugacy of the tower, which is also invariant under
the rescaling:

\begin{prop}\label{prop:17fa,1}
There is a quasi-conformal homeomorphism
$\phi$ of the plane, symmetric w.r.t. the real axis,
and normalized so that $\phi(0)=0, \phi(1)=1, \phi(\infty)=\infty$,
which  conjugates every $H_n$ with $\hat H_n$:
$\phi\circ H_n=\hat H_n\circ \phi$ whenever both sides
are defined.
Moreover, $\phi(z)=\hat \tau \phi(z/\tau)$ for any $z\in {\C}$.
\end{prop}

The conjugacy $\phi$ is easily constructed based on
Proposition~\ref{prop:11fp,1}. Denote $\phi^n(z) = \hat\tau^n
\phi_0(\tau^{-n}z)$. For every $n$, we have 
\[ \phi^n H_n(z) = \hat\tau^n \phi_0 (\tau^{-n}\tau^n H(\tau^{-n} z)) = 
\hat\tau^n \hat{H}(\phi_0(\tau^{-n} z)) = \hat{H}_n(\phi^n(z)) \]
and so $\phi^n$ conjugates $H_n$ to $\hat{H}_n$. Since $H_{n-1} =
H_n^{|\aleph|}$, then $\phi^n$ also conjugates $H_i$ to $\hat{H}_i$
for $i=0,\cdots,n$. 

Using the compactness of the family $\phi^n$, we pick a limit point
$\phi$ which conjugates the whole towers. What will require a check,
however, is the invariance of $\phi$ under the rescaling. 

\paragraph{Uniqueness of the conjugacy on the Julia set.}
\begin{lem}\label{lem:18fa,1}
Suppose that $H$ belongs to the EW-class with some combinatorics
$\aleph$. Let $\Upsilon$ be a homeomorphism which self-conjugates $H$,
i.e. $\Upsilon(H(z)) = H(\Upsilon(z))$ for every $z\in U_-\cup
U_+$. In addition, $\Upsilon$ is symmetric about the real line and
preserves its orientation. 
Then $\Upsilon(z) = z$ for every $z\in K_{H}$. 
\end{lem}
\begin{proof}
We will consider preimages of $[0,x_0]$ by $H^n$ and refer to them as
{\em edges of order k}. The endpoints of each edge of order $n$ are
preimages of $0$: one of order $n$, one of order $n+1$. Let us prove by
induction that $H$ maps each edge of order at most $n$ onto itself,
fixing the endpoints. The first non-trivial case is $n=1$. The edges
of order $1$ are easy to understand: there are two infinite families of
them, one in $U_+$ and one in $U_-$ both branching from $x_0$. 
We can label them
$(e_+^k)_{k=-\infty}^{k=+\infty}$ and $(e_-^k)_{k=-\infty}^{k=+\infty}$,
respectively. See Figure~\ref{fig:18fa,1}. 

\begin{figure}
\epsfig{figure=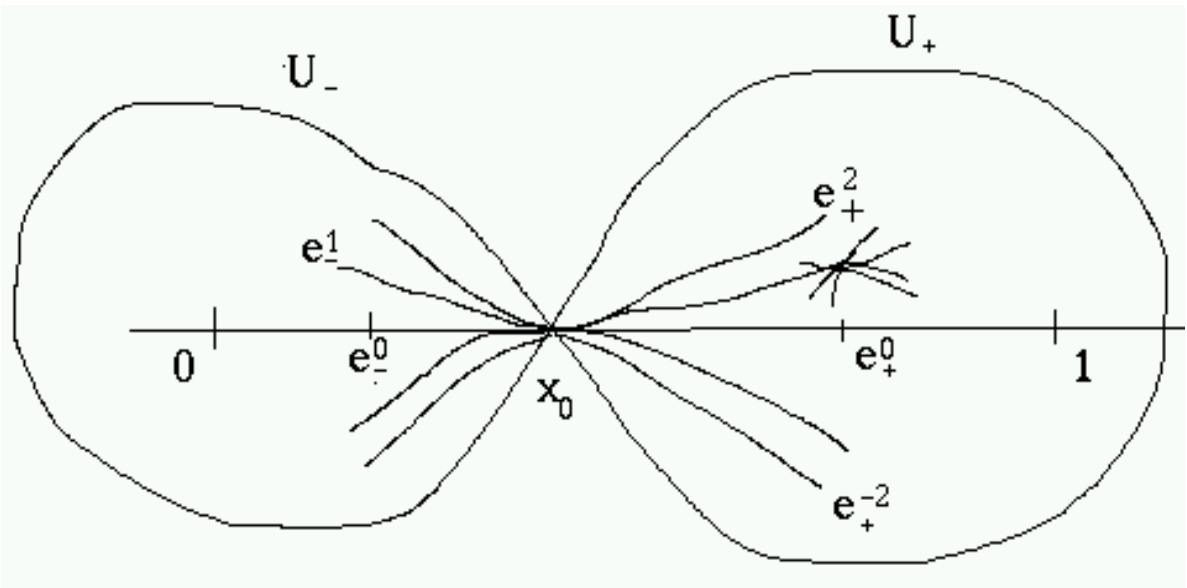}
\caption{Edges of order $1$, with some edges of order $2$ branching
from $e_+^1$.\protect\label{fig:18fa,1}}
\end{figure}

$\Upsilon$ permutes the edges of each family. 
We will focus on the family $e_-^k$ to show that this permutation is
in fact the identity.         
Since $\Upsilon$ preserves the real line with its orientation, we must have
$\Upsilon(e^0_-) = e^0_-$. Then if $e_-^1 = \Upsilon(e_-^{k_1})$ with $k_1 > 1$,
then $\Upsilon(e_-^1)$ would have nowhere to go, since $\Upsilon$ must preserve the
cyclic order of the edges. So $\Upsilon(e_-^1) = e_-^1$ and in this way we
can inductively prove that $\Upsilon(e_-^k) = e_-^k$ for each $k$. 

Now for an inductive step, suppose that $\Upsilon$ fixes all edges of order
$n-1$, $n>1$, but for some edge $e$ of order $n$, $\Upsilon(e) = e' \neq e$. 
One endpoint of $e$ is a preimage of $0$ of order $n$ which also
belongs to an edge of order $n-1$, so it must be fixed by $\Upsilon$. Thus
$e$, $e'$ branch out of the same point $y$ which is the preimage of
$0$ of order $n$. Since $n>1$, a neighborhood of $y$ is mapped by $\Upsilon$
diffeomorphically on to a neighborhood of $\Upsilon(y)$. Then
$\Upsilon(H(e)) = H(\Upsilon(e)) = H(e') \neq H(e)$ which is contrary
to the inductive hypothesis, since $H(e)$ is already an edge of order
$n-1$. 

In particular, it follows that $\Upsilon$ fixes preimages of $0$, but
those are dense in $K_H$ by Theorem~\ref{pre}.  
\end{proof}

Coming back to the proof of Proposition~\ref{prop:17fa,1}, we observe
that for any $m>n$, $\phi^n(z) = \phi^m(z)$ provided that $z\in
K_{H_n}$. Indeed, booth $\phi^m$ and $\phi^n$ conjugate $H_n$ to
$\hat{H}_n$ and so $(\phi^m)^{-1} \circ \phi^n$ provides a
self-conjugacy of $H_n$ and Lemma~\ref{lem:18fa,1} becomes
applicable. 

Now if $\phi = \lim_{k\rightarrow\infty} \phi^{n_k}$, then 
$\phi'(z) = \hat\tau \phi(\tau^{-1} z)$ is the limit of the sequence 
$\phi^{n_k+1}$. For $z\in K_{H_m}$ and any $m$, the values of both
sequences at $z$ stabilize. Hence, $\phi(z) = \phi'(z)$ for any $z\in
\bigcup_{m=0}^{\infty} K_{H_m}$ but this set is dense in $\CC$ by
Proposition~\ref{prop:11fa,3}. So, $\phi=\phi'$ and
Proposition~\ref{prop:17fa,1} has been demonstrated.

\subsection{Invariant line-fields}
We will identify {\em measurable line-fields} with differentials in
the form $\nu(z) \frac{d\overline{z}}{dz}$ where $\nu$ is a measurable
function with values on the unit circle or at the origin. 
A line-field is considered
{\em holomorphic} at $z_0$ if for some holomorphic function $\psi$
defined on a neighborhood of $z_0$, we have $\nu(z) =
c \frac{\overline{\psi'(z)}}{\psi'(z)}$ for some constant $c$.  

By a standard reasoning, Proposition~\ref{prop:17fa,1} gives us a
measurable line-field $\mu(z)\frac{d{\overline z}}{dz}$ which
is invariant under the action of $H_n^*$ for any $n$ as well as under
rescaling: $\mu(\tau z) = \mu(z)$.  
  
We will proceed to show that $\mu$ must be trivial, i.e. $0$ almost
everywhere. This will be attained by a typical approach: showing first
that $\mu$ cannot be non-trivial and holomorphic at any $z_0$ for
dynamical reasons, and on the contrary, that it must be  
holomorphic at some point for analytic reasons and because of
expansion. 

\paragraph{Absence of line-fields holomorphic on an open set.}
\begin{lem}\label{linehol}
The line-field $\mu$ cannot be both holomorphic and non-trivial on any open set.
\end{lem}
\begin{proof} Let $\mu$ be holomorphic in a neighborhood $W$.
Since $\mu$ is invariant under $z\mapsto z/\tau$ and since
$\cup_{n\ge 0} \tau^n K_H$ is dense in the plane,
one can assume that $W$ is a neighborhood of
a point $b$ of $K_H$. Moreover, since $b$ is approximated by
preimages of $x_0$, one can further assume 
that $W$ is a neighborhood of $a$, such that $H^n(a)=x_0$,
for some $n\ge 0$, and (shrinking $W$)
that $H^n$ is univalent on $W$. Apply $H^n$ and see that $\mu$
is holomorphic in a neighborhood $W'$ of $x_0$.
Applying $H$ one more time to $W'\cap U$, one sees
that $\mu$ is holomorphic in a neighborhood
of every point of a punctured disk
$D(0,r)\setminus \{0\}$. Now apply the rescalings
$z\mapsto \tau^n z$, $n=0,1,...$. Hence, 
$\mu$ is holomorphic everywhere except for $0$.
In particular, $\mu$ is holomorphic
around $1=H(0)$. Since $H$ is univalent around $0$,
then $\mu$ is actually holomorphic in the whole disc 
$D(0,r)$. Then $\mu$ cannot be holomorphic around
$H^{-1}(0)=x_0$, a contradiction.
\end{proof}

\paragraph{Construction of holomorphic line-fields.}
Our goal is to prove the following: 
\begin{prop}\label{linetriv}
Suppose that $H$ is a function from the EW-class which fixes an invariant
line-field $\mu(z) \frac{d\overline{z}}{dz}$, which is additionally
invariant under rescaling: $\mu(\tau z) = \mu(z)$. Then the line-field is
holomorphic at some point. Additionally, it is non-trivial in a
neighborhood of the same point unless $\mu(z)$ vanishes almost everywhere. 
\end{prop}

Construction of holomorphic line-fields is based on the following
analytic idea. 

\begin{lem}\label{lem:1gp,1}
Consider  a line-field $\nu_0 \frac{d{\overline
z}}{dz}$ defined on a neighborhood of some point $z_0$ which also is a
Lebesgue (density) 
point for $\nu_0$. Consider a sequence of univalent functions
$\psi_n$ defined on some disk 
$D(z_1,\eta_1)$ chosen so that for every $n$ and a
fixed $\rho<1$ the set $\psi_n(D(z_1,\rho\eta_1))$ covers
$z_0$. In addition, let $\lim_{n\rightarrow\infty} \psi'_n(z_1) = 0$.
Define 
\[ \mu_n(z) \frac{d\overline{z}}{dz}  = \psi_n^*(\nu_0(w))
\frac{d\overline{w}}{dw} \]

Then for some subsequence $n_k$ and a univalent mapping $\psi$ defined
on $D(z_1, \eta_1)$, $\mu_{n_k}(z)$ tend to $\nu_0(z_0)
\frac{\overline{\psi'(z)}}{\psi'(z)}$ on a neighborhood of $z_1$.  
\end{lem}
\begin{proof}
Let us normalize the objects by setting $\hat{\psi}_n :=
|\psi'_n(z_1)|^{-1} \psi_n$ and $\hat{\nu}_n(w) = \nu_0(|\psi'_n(z_1)|
w)$. By bounded distortion, $\hat{\psi}_n(D(z_1,\rho\eta_1))$ contains
some $D(z_0, r_1)$ and is contained in $D(z_0, r_2)$ with $0<r_1 <
r_2$ independent of $n$. By choosing a subsequence, and taking into
account compactness of normalized univalent functions and the fact
that $z_0$ was a Lebesgue point of $\nu_0$, we can assume that
$\hat{\psi}_n$ converge to a univalent function $\psi$ and $\hat{\nu}_n$
converge to a constant line-field $\nu_0(z_0)
\frac{d\overline{w}}{dw}$ almost everywhere. Since 
\[ \mu_n(z) \frac{d\overline{z}}{dz} = \hat{\psi}_n^* (\hat{\nu}_n(w)
\frac{d\overline{w}}{dw} ) \]
for all $n$, we get 
\[ \mu_n(z) frac{d\overline{z}}{dz} \rightarrow \psi^* (\nu_0(z_0)
\frac{d\overline{w}}{dw}) \]
for $z\in D(z_1, \eta_1\rho)$ which concludes the proof of the Lemma. 
\end{proof}

Start with a Lebesgue point $z_0$ of $\mu$. If the field is
non-trivial, without loss of generality $\mu(z_0) \neq 0$. Also, we
can pick $z_0$ so that it is never mapped on the real line and we can
use Proposition~\ref{prop:11fa,1}. 

We then proceed depending on which case occurs in
Proposition~\ref{prop:11fa,1}. In the first case, we choose a point
$Z$ to be an accumulation point of $\tau^{-m_n} z_{n}$. Without loss
of generality, we suppose that $\tau^{-m_n} z_n \rightarrow Z$. The distance
from $Z$ to $\omega_{\infty}$ is positive and we can denote it by
$2\eta_1$. Then, for any $n$ we can find an inverse branch $\zeta_n$ of the
tower iterate $\chi_n$ mapping $z_0$ to $z_{m_n}$ defined on
$D(z_{m_n},\tau^{m_n}\eta)$. One easily checks that functions
$\psi_n(z) = \zeta_n(\tau^{m_n} z)$ defined on $D(Z,\eta)$ satisfy the
hypotheses of Lemma~\ref{lem:1gp,1}. In particular, their derivatives
go to $0$ because $D_{\rho_{\infty}}\chi_n(z_0)$ go to $\infty$ by
Lemma~\ref{lem:5gp,2}.   

To consider the second case of Proposition~\ref{prop:11fa,1}, fix
attention on some $n$. The first observation is that without loss of
generality $|H_{m_n}(z_n)| < R' \tau^{m_n}$ with some $R' < R$
independent of $n$. Indeed, all points on the circle $C(0,
\tau^{m_n} R)$ are in distance $\eta \tau^{m_n}$ from
$\omega_{\infty}$ for some $\eta$ positive. So if this additional
property fails for infinitely many $n$, we can reduce the situation to
the first case already considered. 

Now the key observation is that for every $n$ the point $z_n$ has a
simply connected neighborhood $Y_n$, a point $y_n \in Y_n$ such that
the distance in the hyperbolic metric of $Y_n$ from $z_n$ to $y_n$ is
bounded independently of $n$. Finally, $Y_n$ is  mapped univalently 
by $H_{m_n}$ so that for some integer 
$p_n$ and $\eta>0$ which is independent of
$n$ the image covers $\tau^{p_n}(D(i,\eta))$ with $H_{m_n}(y_n) = \tau^{p_n}
i$. To choose such $Y_n$ and $y_n$, uniformize the component of
$\tau^{m_n}(U_- \cup U_+)$ which contains $z_{m_n}$ by the map 
$\Psi(z) = \log H_{m_n} (z)$ where the branch of the $\log$ is
chosen to make the mapping symmetric about the real axis. $\Psi$ maps
onto the region $\{ \Re w < m_n \log \tau + \log R\}$ and $\Re
\Psi(z_n) < m_n \log\tau + \log R'$. In addition, $|\Im \Psi(z_n)| >
\pi/2$ as the consequence of $\tau^{-m_n} z_n \notin U_{-,c}\cup
U_{+,c}$. Then $Y_n$ can be conveniently chosen in the
$\Psi$-coordinate as a rectangle of uniformly bounded size. 

Once $y_n, Y_n, p_n$ were chosen, we easily conclude the proof. Let
$R_n: D(0,1)\to Y_n$ 
be Riemann maps of regions $Y_n$ with $R_n(0) = y_n$. Then we can
set $\psi_n = (\chi_n)^{-1} \circ R_n$ where $\chi_n$ are maps
specified in Proposition~\ref{prop:11fa,1}. Maps $\psi_n$ satisfy the
conditions of Lemma~\ref{lem:1gp,1}. In particular, 
$|R^{-1}_n(z_n)|$ is bounded independently of $n$ as a consequence of the
construction of $Y_n$.

From this and Proposition~\ref{prop:11fa,1},
the derivatives of $\psi_n$ at $R_n^{-1}(z_n)$ go to $0$, and then the
same can be said of $\psi'_n(0)$ by bounded distortion. So, by passing
to a subsequence, we get that $R_n^* (\mu(z)
\frac{d\overline{z}}{dz})$ tend a.e. to a holomorphic line-field $\nu
\frac{d\overline{w}}{dw}$ on a neighborhood of $0$. 

To finish the proof, we ignore the fact that a subsequence has been
chosen and consider mappings $T_n := \tau^{-p_n} H_{m_n} \circ R_n$
defined on the unit disk. We have $T_n^* (\mu(z)
\frac{d\overline{z}}{dz}) = R_n^* (\mu(z) \frac{d\overline{z}}{dz})$
for every $n$. Maps $T_n$ are all univalent and have been normalized
so that $T_n(0) = i$ and the image of $D(0,1)$ under $T_n$ contains
$D(i, \eta)$ for a fixed $\eta>0$, but avoids $0$. Then $T_n$ is a
compact family of univalent maps and has a univalent limit $T$. Then
it develops that $\mu$ in a neighborhood of $i$ is the image under $T$
of the holomorphic line-field $\nu$ from a neighborhood of $0$, hence
is holomorphic.   

\paragraph{Proof of Theorem~\ref{theo:10fa,2}.}
From Lemma~\ref{linehol} and Proposition~\ref{linetriv} we conclude
that any measurable line-field invariant under the tower of a
EW-mapping and under the rescaling by $\tau$ must be trivial. But as
soon as the conjugacy $\phi$ constructed in
Proposition~\ref{prop:11fa,1} is non-holomorphic, it gives rise to a
non-trivial line field with those properties. Hence, the conjugacy
between any two EW-maps with the same combinatorial pattern must be
holomorphic, and under our normalizations that means the identity.  

This proves Theorem~\ref{theo:10fa,2} which was the last missing link
in the proof of our results.

\section{The Straightening Theorem for EW-maps}
We prove here
\begin{theo}\label{ss}
For every map $H:U_-\cup U_+\to V$ of the EW-class 
there exists a map of the form
$f(z)=\exp(-c(z-a)^{-2})$ with some real $a,c>0$, such that $H$ and
$f$ are hybrid equivalent, i.e.
there exists a quasi-conformal homeomorphism of the
plane $h$, such that
$$h\circ H=f\circ h$$
on $U_-\cup U_+$ and $\partial h/\partial \bar z=0$ a.e. on the 
filled-in Julia set of $H$. 
\end{theo} 
We will see below that $h$ maps the filled-in Julia set 
$K_H$ of $H$ onto the Julia set $J_f$
of $f$.
\begin{proof}
Remind that $V=D(0,R)\setminus \{0\}$.
Making a linear change of variable, one can assume that
$R<1$.
Let us choose real $m>0$, $0<n<R$, as follows.
Consider the map $p(z)=\exp(-m(z-n)^{-2})$, and the set
$\Omega=p^{-1}(V)$.
Then $m,n$ are chosen so that $0\in \Omega$
and $\overline \Omega\subset V\cup \{0\}$.
As in the proof of Proposition~\ref{prop:11fp,2}, one can further choose 
a quasi-conformal homeomorphism $\varphi$ of the plane, such that
$\varphi: V\setminus U_-\cup U_+\to V\setminus \Omega$
is one-to-one, and, most important,
$\varphi(z)=z$ off $V$, and
$\varphi\circ H=p\circ \varphi$ on the boundary
of $U_-\cup U_+$. Also, $\varphi$ is symmetric w.r.t. the real axis.
Since $1\notin V$, we have $\varphi(1)=1$, 
also $\varphi(\infty)=\infty$, and one can assume that
$\varphi(0)=0$.
Now define an extension of $H$ to a map
$\tilde H:{\bf C}\setminus \{x_0\}\to {\bf C}\setminus \{0\}$
as follows:
$\tilde H=H$ on $U_-\cup U_+$, and 
$\tilde H=\varphi^{-1}\circ p\circ \varphi$ on ${\bf C}\setminus(U_-\cup U_+)$.

{\bf Fact 1.} 
{\it Observe that since $\Omega$ is the full preimage of $V$ by $p$,
$p(z)\in {\bf C}\setminus V$ iff 
$z\in {\bf C}\setminus \Omega$.}

Define a complex structure $\sigma$ a.e. on the plane as follows.
Let $\sigma_0$ be the standard one.
Then
$\sigma=\varphi^{*}(\sigma_0)$ on ${\bf C}\setminus \overline{U_-\cup U_+}$;
$\sigma=(H^n)^{*}(\sigma)$ on $H^{-n}(V\setminus \overline{U_-\cup U_+})$,
$n=0,1,2,...$; $\sigma=\sigma_0$ on the rest.
Note that $\sigma=\sigma_0$ off $V$.

As it follows from Fact 1 and since $H$ is holomorphic, we get
 
{\bf Fact 2.} {\it $\sigma$ is correctly defined,
$\tilde H$-invariant, and $||\sigma||_{\infty}<1$.}

Let $h$ be a quasi-conformal homeomorphism of the plane, such that
$h_*(\sigma)=\sigma_0$, $h(0)=0, h(1)=1, h(\infty)=\infty$.
Also, $h$ is symmetric w.r.t. the real axis, because
$\sigma$ is symmetric.
Denote $a=h(x_0)$.
Define $f:{\bf C}\setminus \{a\}\to {\bf C}\setminus \{0\}$ by
$f=h\circ \tilde H\circ h^{-1}$. Then $f$ is holomorphic
because $f_*(\sigma_0)=\sigma_0$.
We need to show that $f(z)=\exp(-c(z-a)^{-2})$, for some real $c>0$.
To this end, notice first that
from the definition of $f$ it follows that
there exists $\lim_{z\to \infty}f(z)=h(1)=1$, and that
$f(z)\not=0$ for every $z\in {\bf C}\setminus \{a\}$.
Hence, the function $\tilde f(z):=1/f(a+1/z)$
is entire. Besides, $\tilde f(z)\not=0$ for any $z$.
Thus there exists another entire function $u$,
such that $\tilde f=\exp(u)$, therefore,
$$f(z)=\exp(-u(1/(z-a))).$$
Let us study singular points of $u^{-1}$
using the formula
$$u^{-1}(w)=[h\circ \tilde H^{-1}\circ h^{-1}(\exp(-w))-a]^{-1}.$$
Since $w\in {\bf C}$, $\exp(-w)\not=0$, hence, $h^{-1}(\exp(-w))\not=0$.
If $\exp(-w_0)\not=1$, then $w_0$ is not a singular point of $u^{-1}$.
If $\exp(-w_0)=1$ but $\tilde H^{-1}\circ h^{-1}(\exp(-w_0))\not=\infty$,
then again $w_0$ is not a singular point.
At last, if $\tilde H^{-1}\circ h^{-1}(\exp(-w_0))=\infty$, then
$w_0$ is a singular point, because then, for $w$ close to $w_0$,
there are two different preimages 
$\tilde H^{-1}\circ h^{-1}(\exp(-w))$ close to $\infty$,
which give two different preimages $u^{-1}(w)$ close to zero. 
Hence, $w_0=0$ is a singular point of $u^{-1}$.
Now, if $w_0=2\pi i k$, $k\in {\bf Z}\setminus \{0\}$, then,
from the symmetry w.r.t. the real axis
and from the continuation along a path $\gamma$ joining $0$ and 
$w_0$, we see using the formula for $u^{-1}$, that
the path $\tilde H^{-1}\circ h^{-1}(\exp(-\gamma))$ is not closed
and starts at $\infty$, hence 
$\tilde H^{-1}\circ h^{-1}(\exp(-w_0))\not=\infty$. Therefore,  
the only singular point of $u^{-1}$ is zero, with the square-root singularity
at this point, and, moreover, $u^{-1}(0)=0$.
Thus, $u(z)=c z^2$, and we are done.
\end{proof}

Now we can make use of the theory of~\cite{bdh},
~\cite{el},~\cite{bkl}
to describe some basic features of the Julia set of $f$.
Remind that the Fatou set $F_f$ is defined as the largest
open set in which all $f^n$ are defined, holomorphic and form a normal family,
and the Julia set $J_f$ is the complement ${\bf \hat{C}}\setminus F_f$.
\begin{prop}\label{j}
Let $V_f=h(V)$ and $U_f=f^{-1}(V_f)=h(U_-\cup U_+)$.
Then:
 
(a) the preimages of the point $a$ are dense in $J_f$, 

(b) $J_f$ is the closure of the set of such $z$ which never leave
$U_f$ under the iterates. 

(c) $J_f$ is connected. 

(d) the Fatou set $F_f={\bf \hat{C}}\setminus J_f$
consists of one component, which is the basin of attraction
of an attractive (real) fixed point of $f$.
Finally, $F_f$ is simply-connected on the sphere.
\end{prop}
\begin{proof}
(a)-(d) follow from a series of 
observations.

(1). The set $E=E(f)$ of singularities of
$f$ consists of one point $a$.
Hence, if $E_n=\cup_{j=0}^{n-1} f^{-j}(E)=\cup_{j=0}^{n-1} f^{-j}(a)$,
then, by~\cite{bdh}, $J_f=\overline{\cup_{n=0}^{\infty} E_n}$.
This proves (a). 

(2). The set $C(f)$ of singular values of $f^{-1}$ consists of the
point $f(\infty)=1$. 

(3). Hence, $f$ belongs to the class
{\bf MSR} defined in~\cite{bdh}.
In particular~\cite{bdh}, 
$f$ has no Baker domains as well as
wandering domains.

(4). We have $b_f:=h(b_0)>a$, and $b_f$ is a repelling
fixed point of $f$. Also, $f$ is strictly increasing on $(a, +\infty)$
and $f(\infty)=1$. Therefore, there exists an attracting fixed point
$z_0$ of $f$, $b_f<z_0<1$.
The iterates of the singular value tend to this fixed point
$z_0$. Hence, for every component $W$ of $F_f$, 
an iterate of $W$ is the immediate basin of attraction $W_0$
of $z_0$.

(5). Since $E_n\subset U_f$ for all $n$,
then, by (1), the domain ${\bf \hat{C}}\setminus \overline{U_f}$ is disjoint 
with $J_f$. It also contains $z_0$. Hence,
${\bf \hat{C}}\setminus \overline{U_f}\subset W_0$.
On the other hand, $z_0\in {\bf C}\setminus V_f$ and
$f^{-1}({\bf C}\setminus V_f)={\bf C}\setminus U_f$, hence,
$f^{-1}(z_0)\subset {\bf C}\setminus U_f\subset W_0$.
Therefore, $W_0$ is completely invariant,
and $F_f=W_0=\cup_{n=0}^\infty f^{-n}({\bf \hat C}\setminus V_f)$.

(6). By (5), $z\in F_f\cap U_f$ iff an iterate of $z$ hits $V_f\setminus U_f$.
Therefore, we have proved that
$J_f=\cap_{n\ge 0} f^{-n}(\overline{U_f})\cup \cup_{n\ge 0}f^{-n}(a)$.
In particular, $J_f$ is connected, and $F_f$ is 
simply-connected. 
\end{proof}
As a corollary, we get a new (indirect) proof of Theorem~\ref{pre}:
\begin{coro}\label{new}
$K_H=h^{-1}(J_f)$, it has no interior, and the preimages
of $x_0$ are dense in $K_H$. 
\end{coro}



Vice versa, one can also gain an information about the dynamics of $f$
from what we know already about the maps $H_\aleph$.
For example, we obtain from Proposition~\ref{prop:11fa,3} that
the union of rescaled (around zero)
Julia sets of $f$ is dense in the plane.
Another information concerns the map $f$ on the real line; let's extend it
to the point $a$ continuously. Then $f:{\bf R}\to {\bf R}$
is a unimodal $C^\infty$ map with the flat critical point
at $a$. Since $f$ and $H_\aleph$ are quasi-conformally conjugate,
then the $\omega$-limit set of the critical point $a$ 
under the dynamics of $f:{\bf R}\to {\bf R}$
has no bounded geometry.

\end{document}